\documentclass{amsart}
\usepackage{amsmath,amsfonts}
\usepackage{amscd,amsthm}
\usepackage{graphicx}
\usepackage{subfigure}

\theoremstyle{plain} 
\newtheorem{theorem}{Theorem}[section]

\newtheorem{proposition}[theorem]{Proposition}

\newtheorem{remark}[theorem]{Remark}

\newcommand{\R}{\mathop{\mathbb{R}}}
\newcommand{\T}{\mathop{\mathbb{T}}}
\newcommand{\Q}{\mathop{\mathbb{Q}}}
\newcommand{\Z}{\mathop{\mathbb{Z}}}
\newcommand{\N}{\mathop{\mathbb{N}}}
\newcommand{\C}{\mathop{\mathbb{C}}}
\newcommand{\B}{\mathop{\mathbf{B}}}
\newcommand{\li}[1]{\mathop{\mathit{Li_2}\left(#1\right)}}
\newcommand{\D}{\mathop{\mathbb{D}}}
\renewcommand{\H}{\mathop{\mathbb{H}^+}}

\numberwithin{equation}{section}

\begin{document}

\title[A numerical test of the Marmi--Moussa--Yoccoz conjecture.]
{The $1/2$--Complex Bruno function and the Yoccoz function. 
A numerical study of the Marmi--Moussa--Yoccoz Conjecture.}

\author{Timoteo Carletti}

\date{\today}

\address[Timoteo Carletti]{Scuola Normale Superiore, piazza dei Cavalieri 7, 56126  Pisa, Italy} 

\email[Timoteo Carletti]{t.carletti@sns.it}



\begin{abstract}
We study the $1/2$--Complex Bruno function and we produce an algorithm
to evaluate it numerically, giving a characterization of the
monoid $\hat{\mathcal{M}}=\mathcal{M}_T\cup \mathcal{M}_S$.
We use this algorithm to test the Marmi--Moussa--Yoccoz Conjecture about the 
H\"older continuity of the function $z\mapsto -i\B(z)+
\log U\!\left(e^{2\pi i z}\right)$ on $\{ z\in \C: \Im z \geq 0 \}$, where $\B$ is the $1/2$--complex 
Bruno function and $U$ is the Yoccoz function. We give a positive
answer to an explicit question of S. Marmi et al~\cite{MMYc}.
\end{abstract}

\maketitle

\section{Introduction}
\label{sec:intro}

The {\em real Bruno functions} are arithmetical functions 
$B_{\alpha}:\R \setminus \Q \rightarrow \R_{+} \cup \{ +\infty \}$, 
$\alpha\in[1/2,1]$ which characterize numbers by their rate of 
approximation by rationals. They have been introduced 
by J.--C. Yoccoz~\cite{Yoccoz1} (cases $\alpha=1/2$ and $\alpha=1$)
and then 
studied in a more general
context in~\cite{MMYr}.

For their relationship with arithmetical properties of real numbers, 
Bruno's functions enter in a huge number of dynamical
system problems involving small divisors, for instance in the 
problem of the {\em stability
of a fixed point} of a holomorphic diffeomorphism of a complex variable
(the so called {\em Schr\"oder--Siegel} problem)~\cite{Yoccoz1},
in the {\em Schr\"oder--Siegel} problem in the {\em Gevrey} setting in
one complex variable~\cite{CarlettiMarmi} and in some {\em local conjugacy}
 problems: {\em Semistandard map}~\cite{Marmi,Davie}, 
{\em analytic circle diffeomorphisms}~\cite{Yoccoz2} and some
 {\em analytic area--preserving annulus map} including
the {\em Standard map} and some of its generalizations~\cite{BerrettiGentile}.

Let us now concentrate on the $1/2$--Bruno function~\footnote{From~\cite{MMYr} 
we know that the difference of any two Bruno's functions is in 
$L^{\infty}(\R)$.}.
 $B_{1/2}$ is $\Z$--periodic, even (for this reason it
is also 
called
 {\em even Bruno's function}) and verifies the {\em functional equation}:
\begin{equation}
\label{eq:funceqBr}
B_{1/2}\left(x\right)=-\log x +xB_{1/2}\left(x^{-1}\right)\quad x\in
(0,1/2) \, .
\end{equation}
The set $\mathcal{B}=\{ x\in \R: B_{1/2}\left( x\right)<+\infty\}$ 
is called the set {\em of Bruno's numbers}: by~\eqref{eq:funceqBr} it
follows that $\mathcal{B}$ is invariant under the action of the
modular group $GL\left(2,\Z\right)=\Big\{\left(\begin{smallmatrix}a &b \\ c &d
\end{smallmatrix}\right): a,b,c,d \in \Z, ad-bc =\pm 1 \Big\}$. The Bruno function
can be extended to rational numbers by setting $B_{1/2}(x)=+\infty$ when $x\in \Q$.

 Using the {\em continued fraction} algorithm one can
 solve~\eqref{eq:funceqBr} to obtain:
\begin{equation}
\label{eq:defBr}
B_{1/2}\left(x\right)=\sum_{k\geq 0}\beta_{k-1}\left(x\right) \log
x_k^{-1} \, , 
\end{equation}
where $x_0=x$, $x_k=A_{1/2}\left( x_{k-1} \right)$, $\beta_{-1}=1$,
$\beta_k=\prod_{j=0}^k x_j$ 
and $A_{1/2}$ is the {\em nearest integer Continued Fraction map}. In
\S~\ref{ssec:contfrac} we will give a brief account of useful facts
concerning continued fractions.

In~\cite{MMYc} the {\em complex
  Bruno function} has been introduced~\footnote{ 
Following the notation introduced for the real Bruno functions, we
should call this  
complex extension the {\em $1$--complex Bruno function}. In fact we
will see at the end 
of section~\ref{sec:cmplxBruno}, that it is constructed ``following'' the Gauss
continued fraction algorithm. In this way we could also distinguish it
  from the  
{\em $1/2$--complex Bruno function} that we will introduce in
  section~\ref{sec:cmplxBruno} 
``following'' the nearest integer continued fraction algorithm.}, more
precisely authors defined an analytic map $\B:\H \rightarrow \H$, 
 where $\H$ is the {\em upper Poincar\'e half plane},
$\Z$--periodic, which verifies a functional equation similar to the one for the
$1$--Bruno function. The boundary  behavior of $\B$ is given by (see Theorem 5.19
and \S 5.2.9 of~\cite{MMYc}):
%
\begin{enumerate}
\item let $H>0$, then the imaginary part of $\B(z)$ tends to $B_{1/2}(\Re z)$ when 
$\Im z \rightarrow 0$ and $z\in \{ \zeta \in \H : \Im \zeta \geq | \Re \zeta |^H\}$,
whenever $\Re z \in \mathcal{B}$;
\item $\Re \B(z)$ is bounded on $\H$, its trace 
is continuous at irrational points and it has a jump of
$\pi/q$ for $\Re z=p/q \in \Q$.
\end{enumerate}
In section~\ref{sec:cmplxBruno} we introduce an explicit formula for
the $1/2$--complex Bruno function which corrects a small error in 
\S A.4.4 page 836, and gives more details than Appendix
A.4 of~\cite{MMYc}. We will also give an algorithm to compute it numerically. 

\subsection{The Yoccoz function}
\label{ssec:yoccozfunc}

We already observed that the function $B_{1/2}$ is related to the
stability problem of a fixed point of an analytic diffeomorphism of
$\C$, in the rest of this section we will show this relation by
describing the Yoccoz result (\cite{Yoccoz1}, Chapter II). Let $\lambda
\in \C^*$ and 
let us consider the {\em quadratic polynomial} $P_{\lambda}(z)=\lambda
z(1-z)$. The origin is a fixed point and we are interested in studying
its stability. If $|\lambda|<1$ (hyperbolic case) then it follows from
the results 
of Poincar\'e and Koenigs that the
origin is stable, whereas if $\lambda=e^{2\pi i p/q}$ (parabolic case)
the origin is not stable.

Let now consider $\lambda \in \D^*$ and let $H_{\lambda}(z)$ be the
conformal map 
which locally linearizes $P_{\lambda}$ (its existence is guaranteed by
the Poincar\'e--Koenigs results):
\begin{equation}
  P_{\lambda}\circ H_{\lambda}= H_{\lambda}\circ R_{\lambda}\, ,
\label{eq:siegel}
\end{equation}
where $R_{\lambda}(z)=\lambda z$, and let us denote by $r_2(\lambda)$ the
radius of convergence of $H_{\lambda}$. 

One can prove that $H_{\lambda}$ can be analytically continued to a
larger set, the {\em basin of attraction of $0$}: $\{ z \in
\C:P_{\lambda}^{\circ n}(z)\rightarrow 0, n\rightarrow +\infty \}$,
but not to the whole of $\C$, and it has a unique singular point on
its circle of convergence $\D_{r_2(\lambda)}$, which will be denoted
by $U\!(\lambda) \in \C$. The function  $U:\D^* \rightarrow \C$ is 
called the {\em  Yoccoz function}.

Yoccoz proved that $U$ has an analytic bounded extension to $\D$ and
moreover it can be obtained as limit of polynomials
$U_n(\lambda)=\lambda^{-n}P^{\circ n}_{\lambda}\left(z_{crit}\right)$,
uniformly over compact subsets of $\D$, where $z_{crit}=1/2$ is the
critical point of the quadratic polynomial. Being this extension not
identically zero, by a classical result of Fatou, the
Yoccoz function has radial limits almost everywhere, and the set
$\lambda_0 \in S^1$ for which $\limsup_{\lambda \rightarrow \lambda_0}
U\!(\lambda) = 0$ has zero measure. This means that the quadratic
polynomial is linearizable ($r_2(\lambda_0)=|U\!(\lambda_0)| >0$)
for a full measure set of $\lambda_0 \in
S^1$, but the proof doesn't give any information on this set.

When $|\lambda|=1$ and $\lambda$ is not a {\em
  root of the unity}, assuming $\lambda=e^{2\pi i \omega}$, for
some irrational $|\omega |<1/2$, Yoccoz proved~(\cite{Yoccoz1},
Theorem 1.8 Chapter II) that
$P_{\lambda}(z)$ is linearizable if and only if $\omega \in
\mathcal{B}_{1/2}$, moreover there exists a constant
$C_1$, and for all $\epsilon >0$ a constant $C(\epsilon)$ such 
that for all $\omega \in \mathcal{B}_{1/2}$:
\begin{equation*}
  C_1 \leq \log r_2(e^{2\pi i \omega}) +B_{1/2}(\omega) \leq
  C(\epsilon)+\epsilon B_{1/2}(\omega) \, .
\end{equation*}
We are then interested in studying the function $\omega \mapsto \log
  |U\left(e^{2\pi i 
  \omega}\right)| +B_{1/2}(\omega)$ and some ``natural''
 questions arise (\cite{Yoccoz1} \S 3.2 page 72):
  \begin{quotation}
    (Yoccoz's Conjecture). {\em Is the function $\omega \mapsto
  \log |U\left(e^{2\pi i \omega}\right)| +B_{1/2}(\omega)$ bounded for $\omega\in\R$?}
  \end{quotation}
Motivated by numerical results of~\cite{Marmi} and by some analytic 
properties of the real Bruno function (see the following 
Remark~\ref{rem:12hold} and~\cite{MMYr}) it has been conjectured that:
  \begin{quotation}
    (Marmi--Moussa--Yoccoz's Conjecture). {\em The function, defined on the set
 of Bruno number, $\omega \mapsto \log |U\left(e^{2\pi i \omega}\right)| 
+B_{1/2}(\omega)$, extends to a $1/2$--H\"older continuous function on $\R$.}
  \end{quotation}

\begin{remark}[Why $1/2$--H\"older?]
  \label{rem:12hold}
In~\cite{MMYr} authors proved a ``stability result'' for $B_{1/2}$ (\S
4 page 285). Let us rewrite the functional equation for the 
$1/2$--Bruno function as follow:
\begin{equation*}
\left[B_{1/2}(x)-xB_{1/2}\left(x^{-1}\right)\right]=-\log x \, ,
\end{equation*}
if we add to the r.h.s. a ``regular term'' $f$, say
$\eta$--H\"older continuous, and we call $B_f$ the solution of:
\begin{equation*}
\left[B_{f}(x)-xB_{f}\left(x^{-1}\right)\right]=-\log x +f(x) \, ,
\end{equation*}
then $B_{1/2}-B_f$ is $1/2$--H\"older continuous if $f$ is
at least $1/2$--H\"older. Hence if we prove~\footnote{Transform a function according to
$\psi(x) \mapsto \psi(x)-x\psi(1/x)$ to ``reduce the strength of singularities'' is the main
idea of the {\em Modular Smoothing}. We refer to~\cite{BuricPercivalVivaldi} where authors 
describe the method and apply it to the critical function of the Semistandard Map.} that the function
$\omega \mapsto \left[\log |U\left(e^{2\pi i \omega}\right)|-\omega 
\log \Big\lvert U\left(e^{2\pi i \omega^{-1}}\right)\Big\rvert\right]-\log \omega$
is H\"older continuous with exponent $\eta \geq 1/2$, for $\omega \in [0,1/2]$,
 then the Conjecture holds.
\end{remark}

Very recently X. Buff and A. Cheritat~\cite{BuffCheritat} proved the Yoccoz
conjecture. Whereas the Marmi--Moussa--Yoccoz conjecture is still open.
We will be interested in the following conjecture, equivalent to the
one of Marmi--Moussa--Yoccoz:
  \begin{quotation}
    (Marmi--Moussa--Yoccoz's Conjecture)
{\em The analytic function, defined on the upper Poincar\'e half plane, 
   $z \mapsto \mathcal{H}\left( z\right)=
  \log U \! (e^{2\pi i z}) -i \B(z)$, extends to a
  $1/2$--H\"older continuous function on $\bar{\mathbb{H}}^{+}$.}
  \end{quotation}

The aim of this paper is twofold. First give more insight into the
$1/2$--complex Bruno function and second to make a first step toward
the understanding of the last conjecture. Our numerical results 
allows to conclude that $\mathcal{H}$ is $\eta$--H\"older continuous 
and we obtain an estimate of the H\"older exponent $\eta = 0.498\pm
0.004$. This 
gives us good numerical evidence that the Marmi--Moussa--Yoccoz conjecture should be true.

The paper is organized as follows: in section~\ref{sec:cmplxBruno} we 
introduce the $1/2$--complex Bruno function and some results from
number theory (approximations of rationals by rationals), to obtain
 an algorithm to compute the complex Bruno
function. In section~\ref{sec:yoccozf} we explain how to calculate
 the Yoccoz function and then, after a brief introduction of the
Littlewood--Paley Theory in \S~\ref{sec:lpm}, used to test the H\"older continuity,
 we present our results in section~\ref{sec:numrel}. Appendix~\ref{sec:numcons} collects some
considerations 
related to technical aspects of our numerical test.

{\it Acknowledgments.}
I am grateful to Jacques Laskar for putting at my disposal a large amount
of CPU times on computers
of {\em Astronomie et Syst\`emes Dynamiques} Team at IMCCE Paris.

\section{The $1/2$--Complex Bruno functions}
\label{sec:cmplxBruno}

The aim of this section is to introduce, starting from Appendix A.4 of~\cite{MMYc},
 a complex extension of the $1/2$--real Bruno function and to give an algorithm
to compute it numerically.

Let us consider $f\in L^2([0,1/2])$, extended: $1$--periodic,
 $f(x+1)=f(x)$ for all $x\in \R$, 
and even $f(x)=f(-x)$ for all $x\in [-1/2,0]$, and then let us introduce the 
operator $T$ acting on such $f$ by:
\begin{equation}
\label{eq:teven}
Tf(x)=xf\left( \frac{1}{x}\right) \, ;
\end{equation}
we remark that the functional equation~\eqref{eq:funceqBr}
 can be rewritten as:
\begin{equation}
\left(1-T\right)B_{1/2}(x)=-\log x \quad \forall x\in (0,1/2) \, .
\label{eq:tbfunc}
\end{equation}

Let $(T_m)_{m\geq 2}$ be the operators defined by:
\begin{equation}
  \label{eq:tmop}
  \left(T_mf\right)(x)=\begin{cases}x f\left(\frac{1}{x}-m\right) & x\in \left( 
\frac{1}{m+1/2},\frac{1}{m} \right] \text{branch $m^+$} \\
x f\left(m-\frac{1}{x}\right) & x\in \left( 
\frac{1}{m},\frac{1}{m-1/2} \right] \text{branch $m^-$}\\
0 & \text{otherwise,}
\end{cases}
\end{equation}
then using the periodicity and the evenness of $f$, we 
can rewrite~\eqref{eq:teven} as follows:
\begin{equation}
\label{eq:tevenm}
Tf(x)=\sum_{m\geq 2} \Big\{ xf\left(\frac{1}{x}-m\right)+
xf\left( m+1-\frac{1}{x}\right) \Big \} \, .
\end{equation}

To introduce the $1/2$--complex Bruno function we have to extend~\eqref{eq:tevenm} to 
complex analytic functions; this is done~\cite{MMYc} by considering the 
complex vector space of holomorphic functions in $\bar{\C}\setminus [0,1/2]$,
 vanishing at infinity: $\mathcal{O}^{1}\left(\bar{\C}\setminus [0,1/2]\right)$
 (which is isomorphic to the space of hyperfunctions
with support contained in $[0,1/2]$). So, let $\varphi$ be the Hilbert transform of $f$:
\begin{equation*}
\varphi(z)=\frac{1}{\pi}\int_0^{1/2} \frac{f(x)}{x-z} \, dx \, ,
\end{equation*}
then starting from~\eqref{eq:tevenm} we define the action of $T$ on $\varphi$ as follows:
\begin{equation}
\label{eq:complexT}
T\varphi(z)=\sum_{m\geq 2}L_{g(m)}\left( 1+L_{\sigma}\right) \varphi(z) \, ,
\end{equation}
where $g(m)=\left(\begin{smallmatrix} 0 & 1 \\ 1 & m\end{smallmatrix}\right)$,
 $\sigma=\left(\begin{smallmatrix} -1 & 1 \\ 0 & 1\end{smallmatrix}\right)$ and
$L_{\left(\begin{smallmatrix} a & b \\ c & d\end{smallmatrix}\right)}$ acts
 on $\mathcal{O}^{1}\left(\bar{\C}\setminus [0,1/2]\right)$ by:
\begin{equation}
\label{eq:Laction}
L_{\left(\begin{smallmatrix} a & b \\ c & d\end{smallmatrix}\right)}
  \varphi(z)=(a-cz)\left[ \varphi\left(
    \frac{dz-b}{a-cz}\right)-\varphi\left(
    -\frac{d}{c}\right)\right]-\frac{ad-bc}{c} \varphi^{\prime}\left(
    -\frac{d}{c}\right) \, .
\end{equation}
In the spirit of~\eqref{eq:tbfunc} we want to consider $\left(1-T\right)^{-1}$ 
acting on some $\varphi \in \mathcal{O}^{1}\left(\bar{\C}\setminus
[0,1/2]\right)$, and to 
obtain a $\Z$--periodic, ``even function''~\footnote{Here and in the
  following by even complex function we will mean even w.r.t. $\Re z
  \rightarrow -\Re z$.}, we will consider:
\begin{equation}
\label{eq:complex1}
\sum_{n\in
  \Z}\left[\left(1+L_{\sigma}\right)\left(1-T\right)^{-1}\right]\varphi
  (z-n) 
\, . 
\end{equation}
Let us introduce the operator $\hat{T}$ defined by
$\left(1+L_{\sigma}\right)T=\hat{T}\left(1+L_{\sigma}\right)$, 
then from~\eqref{eq:complexT} and the relation:
$\left(1-T\right)^{-1}=\sum_{r\geq 0}T^r$,
we can expand $\left(1-\hat{T}\right)^{-1}$  in terms of matrices
$g(m)$ and $\sigma$, to obtain a sum of matrices of the form: 
$\epsilon_0 g(m_1) \dots \epsilon_{r-1}g(m_r)$, where $r\geq 1$,
$m_i \geq 2$ and $\epsilon_{i-1} \in \{ 1,\sigma \}$, for $1\leq i \leq r$.

Let us set $\hat{\mathcal{M}}^{(0)}=\{ 1 \}$ and for $r\geq 1$:
\begin{eqnarray}
\hat{\mathcal{M}}^{(r)}=\Big\{ g\in GL(2,\Z) : \exists \epsilon_0,\dots,\epsilon_{r-1} \in \{1,\sigma\},
m_1,\dots,m_r  \geq 2 & : \notag \\  
g= \epsilon_0 g(m_1) \dots \epsilon_{r-1}g(m_r)&\Big\} \, ,
\label{eq:monoidr}
\end{eqnarray}
and finally $\hat{\mathcal{M}}=\cup_{r\geq 0} \hat{\mathcal{M}}^{(r)}$:
the {\em $1/2$--Monoid} (we left to \S~\ref{ssec:12monoid} a more detailed discussion
of this monoid and the reason of its name).

It remains to specify the ``good'' $\varphi \in \mathcal{O}^{1}\left(\bar{\C}\setminus [0,1/2]\right)$
to apply~\eqref{eq:complex1}, to have the wanted properties for $\B$. This is done
by considering the Hilbert transform of the {\em even, $1$--periodic} logarithm
 function defined in $(0,1/2]$, namely:
\begin{eqnarray}
\label{eq:htrasnlog}
\varphi_{1/2}(z)&=&\frac{1}{\pi}\int_0^{1/2}\frac{-\log x}{x-z}\, dx \notag \\
             &=&-\frac{1}{\pi}\li{\frac{1}{2z}}+\frac{1}{\pi}\log 2 \log \left(1-
\frac{1}{2z} \right) \, ,
\end{eqnarray}
where $\li{z}$ is the {\em dilogarithm} function~\cite{Oesterle}: the analytic continuation of
$\sum_{n\geq 1}z^n n^{-2}$, to $\C\setminus [1,+\infty)$. We are now able to define 
the {\em $1/2$--complex Bruno function} to be:
\begin{equation}
\label{eq:defcomplxB}
\B(z)=\sum_{n\in \Z}\left[\sum_{g\in\hat{\mathcal{M}}}L_g\left(1+L_{\sigma}\right)\right] 
\varphi_{1/2}(z-n) \, .
\end{equation}
This formula defines~\footnote{This claim can be obtained by slightly modification
 of the proof given in~\cite{MMYc} for the $1$--Complex Bruno function and we
omit it referring to~\cite{MMYc} for any details.} an holomorphic function,
defined in $\H$, $\Z$--periodic, with an {\em even imaginary part} for
$\Re z\in [0,1/2]$, and then~\footnote{This is a standard result for
 harmonic conjugate functions. For this reason we prefer to speak of
 $1/2$--Complex Bruno's function instead of even Complex Bruno's
 function.} an {\em odd real part}, on the same domain. 

\begin{remark}
\label{rem:fourcoeffB}
The $1/2$--complex Bruno function is $1$--periodic and so we can consider
its Fourier series: $\B\left(z\right)=\sum_{l\in \Z}\hat{b}_l e^{2\pi il z}$.
Introducing the variable $w=e^{2\pi i z}$ the Bruno function is mapped 
into an analytic function, $\Tilde{\B}\left(w\right)$, defined in $\D^*$,
which can be extended by continuity to $\D$. Its Taylor series at
the origin is $\Tilde{\B}\left(w\right)=\sum_{l\in \N}\hat{b}_l w^{l}$,
hence Fourier coefficients of $\B\left(z\right)$ corresponding
to negative modes are all identically zero. Moreover, because of the parity properties of $\Re\B$ and $\Im\B$, its Fourier 
coefficients are all purely imaginary, in fact:
\begin{equation*}
\hat{b}_l=2i\int_0^{1/2} \left[ -\sin \left(2\pi lx \right)\, \Re
  \B\left(x+it\right)+ 
\cos \left(2\pi lx \right)\, \Im \B\left(x+it\right)\right] \, dx \, .
\end{equation*}
\end{remark}

The goal of the next sections will be to express~\eqref{eq:defcomplxB} in terms of
a sum over a class of rational numbers in such a way we could give 
(\S~\ref{ssec:cmplxbruno}) an algorithm to compute it. This will be 
accomplished thanks to a new characterization (\S~\ref{ssec:12monoid}) of the
{\em $1/2$--Monoid} $\hat{\mathcal{M}}$, after having introduced some results
from Number Theory: \S~\ref{ssec:contfrac} Continued Fraction Theory 
 and \S~\ref{ssec:approx} Farey Series.

\subsection{Continued Fraction}
\label{ssec:contfrac}

We consider the so called {\em nearest integer Continued Fraction} 
algorithm~\footnote{In~\cite{MMYr} a one parameter family of continued
fraction developments has been introduced. The nearest integer continued fraction corresponds to
the value $1/2$ of the parameter, so we will also call it $1/2$--continued fraction.}.
 We state here some basic facts we will need in the following and we 
refer to~\cite{HardyWright,MMYr}
 for a more complete discussion. Let $\lvert\lvert x\rvert\rvert=
\min_{p\in\Z} \{ x<1/2+p \}$, then to each $x\in \R$ we associate a 
continued fraction as follows:
\begin{equation}
\label{eq:cfa0}
a_0 =\lvert\lvert x\rvert\rvert \, , \quad x_0 =|x-a_0| \, , \quad 
\varepsilon_0=\begin{cases}+1 \text{ iff } x\geq a_0 \\
                        -1 \text{ otherwise}
           \end{cases} \, ,
\end{equation}
and then inductively for all $n\geq 0$, as long as $x_n \neq 0$:
\begin{equation}
\label{eq:cfan}
a_{n+1} =\lvert\lvert x_n^{-1}\rvert\rvert \, , \quad x_{n+1}
=|x_n^{-1}-a_{n+1}|\equiv A_{1/2}(x_n) \, , \quad  
\varepsilon_{n+1}=\begin{cases}+1 \text{ iff } x_n^{-1}\geq a_{n+1} \\
                        -1 \text{ otherwise}
           \end{cases} \, .
\end{equation}
We will use the standard compact notation to denote the continued fraction 
$x=[(a_0,\varepsilon_0),\dots,(a_n+\varepsilon_nx_n,\varepsilon_n)]$.
 From the definition it follows that $x_n > 2$ and so $a_n \geq 2$.

\begin{remark}[Standard form for finite continued fraction]
\label{rem:conv}
Let $[(a_0,\varepsilon_0),\dots,(a_{\bar{n}},\varepsilon_{\bar{n}})]$ be a
finite continued fraction of length $\bar{n}$. Then, whenever 
$a_{\bar{n}}=2$, we must also have $\varepsilon_{\bar{n}-1}=+1$, namely
$[(a_0,\varepsilon_0),\dots,(a_{\bar{n}-1},-1),(2,+1)]$ represents the
same rational number that
$[(a_0,\varepsilon_0),\dots,(a_{\bar{n}-1}-1,+1),(2,+1)]$. Moreover a finite 
continued fraction cannot contain a couple $(a_l,\varepsilon_l)=(2,-1)$ 
for any $l\leq \bar{n}$.
\end{remark}

We recall, without proof some known results:
\begin{itemize}
\item the continued fraction algorithm stops if and only if $x\in
  \R\setminus\Q$, this correspondence in bijective up to the standard
  convention of Remark~\ref{rem:conv};
\item \label{it:pnqn}For any positive integer $n$ (or smaller than the length of the 
finite continued fraction) the $n^{\text{th}}$ convergent is defined by:
\begin{equation}
\label{eq:pnqn}
\frac{p_n}{q_n}=[(a_0,\varepsilon_0),\dots,(a_n,\varepsilon_n)] \, ,
\end{equation}
one can  prove that $p_n$ and $q_n$ are recursively defined by:
\begin{equation}
\label{eq:pnqn2}
\begin{cases}
p_n = a_n p_{n-1}+\varepsilon_{n-1}p_{n-2} \\
q_n = a_n q_{n-1}+\varepsilon_{n-1}q_{n-2} \, ,
\end{cases}
\end{equation}
starting with $p_{-1}=q_{-2}=1$, $p_{-2}=q_{-1}=0$ and $\varepsilon_{-1}=1$.
\item \label{it:segnopnqn} for all $n$ we have: 
$q_np_{n-1}-p_nq_{n-1}=(-1)^n\varepsilon_0 \dots \varepsilon_{n-1}$.
\end{itemize}

\subsection{The Farey Series}
\label{ssec:approx}

Let $n\in \N^*$, the {\em Farey Series}~\cite{HardyWright} of order
$n$ is the set of irreducible 
fractions in $[0,1]$ whose denominators do not exceed
$n$~\footnote{This is different from the {\em Farey Tree} which is  
still a set of rational numbers in $[0,1]$ which can be constructed
by induction starting with:  
$\hat{\mathcal{F}}_0=\{0,1\}$ and then defining the $i$--th element of
$\hat{\mathcal{F}}_n$, $n\geq 1$, by: 
\begin{equation*}
\frac{\hat{p}^{(n)}_i}{\hat{q}^{(n)}_i}=\frac{\hat{p}^{(n-1)}_{i-1}+
\hat{p}^{(n-1)}_{i}}{\hat{q}^{(n-1)}_{i-1}+\hat{q}^{(n-1)}_{i}}
\, . 
\end{equation*} 
The Farey Tree of order $n$ is clearly larger than the corresponding
 Farey Series and $card \hat{\mathcal{F}}_n = 2^n+1$.}:
\begin{equation}
  \label{eq:fareyseries}
  \mathcal{F}_n=\{ p/q \in [0,1] : (p,q)=1 \text{ and } q\leq n \}.
\end{equation}
The cardinality of $\mathcal{F}_n$ is given by the Euler $\phi(n)$
function and so it is asymptotic to $3n^2/\pi^2$ for $n$ large. 
The Farey Series is characterized by the following two
properties~\cite{HardyWright}: 
\begin{theorem}
\label{thm:1carat}
  Let $n\geq 1$. If $p/q$ and $p^{\prime}/q^{\prime}$ are two successive 
elements of $\mathcal{F}_n$, then: $qp^{\prime}-q^{\prime}p=1$.
\end{theorem}

\begin{theorem}
\label{thm:thm2}
Let $n\geq 1$. If $p^{\prime}/q^{\prime}$, $p/q$ and 
$p^{\prime\prime}/q^{\prime\prime}$ are three successive elements (in
this order) 
 of $\mathcal{F}_n$, then: 
 \begin{equation*}
   \frac{p}{q}=\frac{p^{\prime}+p^{\prime\prime}}{q^{\prime}+q^{\prime\prime}}
   \, . 
 \end{equation*}
\end{theorem}

Using an idea contained in the proof of Theorem~\ref{thm:thm2} given 
in~\cite{HardyWright}, 
 we construct an algorithm (easily implementable on a computer)
 which allows us to carry out for any $n \geq 2$ the Farey
 Series of order $n$. Using Proposition~\ref{prop:farey2} we will give a second
algorithm to compute the Farey series up to any given order $n$,
 using the continued fraction development.

\begin{proposition}[Construction of $\mathcal{F}_n$]
Let $n\geq 2$, then the elements of $\mathcal{F}_n$, 
$(p_i/q_i)_{1\leq i \leq \phi(n)}$, are recursively defined by:
\begin{equation}
\label{eq:constrFn}
\begin{cases}
p_{i+1}=-p_{i-1}+r_ip_i \\
q_{i+1}=-q_{i-1}+r_iq_i \, ,
\end{cases}
\end{equation}
where $r_i=\lfloor (n+q_{i-1})/q_i\rfloor$, starting with $(p_1,q_1)=(0,1)$, 
$(p_2,q_2)=(1,n)$ and $(p_3,q_3)=(1,n-1)$.
\end{proposition}

\proof
Let $p/q \in \mathcal{F}_n$. Because $p$ and $q$ are relatively prime we can 
always solve in $\Z^2$ the linear Diophantine equation $qP-pQ=1$: let
$(P_0,Q_0)$ be a particular 
 solution and let $r$ be the integer such that: $n-q<Q_0+rq\leq n$, namely
$r=\lfloor (n-Q_0)/q\rfloor$.

Let us define $P_r=P_0+rp$ and $Q_r=Q_0+rq$, then the following claims
are trivial: $(P_r,Q_r)$ is again a solution of the linear diophantine
equation, $(P_r,Q_r)=1$ and $0< Q_r\leq n$. So $P_r/Q_r\in\mathcal{F}_n$. Clearly 
$P_r/Q_r > p/q$ and we claim that it is the immediate
successor of $p/q$ in $\mathcal{F}_n$.

To obtain a constructive algorithm we must solve the linear diophantine 
equation, this is achieved by considering the element which precedes $p/q$ 
in $\mathcal{F}_n$: let us denote it by $p^{\prime}/q^{\prime}$. A
particular solution is then 
given by $P_0=-p^{\prime}$, $Q_0=-q^{\prime}$:
 from the previous result the element following $p/q$ is then given by:
$P_r=-p^{\prime}+rp$, $Q_r=-q^{\prime}+rq$, where $r=\lfloor
(n+q^{\prime})/q\rfloor$. 

To finish the algorithm we need two starting elements of
 $\mathcal{F}_n$ apart of $0/1$, 
 but it is easy to realize that the first three elements of $\mathcal{F}_n$ 
are $0/1$, $1/n$ and $1/(n-1)$, whenever $n\geq 2$.
\endproof

We are now able to give a second algorithm to construct the Farey
Series of order $n$. Here is the idea: given an irreducible fraction
$p/q\in(0,1)$, we compute its continued fraction development and then
following two rules: {\em Truncate} and {\em Subtract one},
we obtain two new irreducible fractions in $[0,1]$ which will result to be
the predecessor and the successor of $p/q$ in $\mathcal{F}_n$ with $n=q$.

\begin{proposition}[Construction of $\mathcal{F}_n$, $2^{\text{nd}}$ version]
\label{prop:farey2}
Let $p/q \in (0,1)$ and let
$p^{\prime}/q^{\prime}<p/q<p^{\prime\prime}/q^{\prime\prime}$ 
be three successive elements of $\mathcal{F}_q$. Assume 
$p/q=[(a_0,\varepsilon_0),\dots,(a_{\bar{n}},\varepsilon_{\bar{n}})]$ for
some $\bar{n}\geq 1$ and 
let us define the rational numbers $p_T/q_T$ and $p_S/q_S$ as 
follows~\footnote{If $a_{\bar{n}}=2$, then $\varepsilon_{\bar{n}-1}=+1$
  by remark~\ref{rem:conv}, 
and $p_S/q_S=[(a_0,\varepsilon_0),\dots,(a_{\bar{n}-1}+1,+1)]$.}:
\begin{equation}
\label{eq:ptqt}
\frac{p_T}{q_T}=[(a_0,\varepsilon_0),\dots,(a_{\bar{n}-1},\varepsilon_{\bar{n}-1})] 
\quad \text{(Truncate)} \, ,
\end{equation}
and 
\begin{equation}
\label{eq:pmqm}
\frac{p_S}{q_S}=[(a_0,\varepsilon_0),\dots,(a_{\bar{n}}-1,\varepsilon_{\bar{n}})] 
\quad \text{(Subtract one)} \, .
\end{equation}
Then if $\varepsilon_0 \dots \varepsilon_{\bar{n}-1}=+1$, we have
\begin{equation}
\label{eq:resp1}
\begin{cases}
p_T/q_T = p^{\prime}/q^{\prime} \, , \text{and,} \, \,
p_S/q_S = p^{\prime\prime}/q^{\prime\prime} \quad \text{if $\bar{n}$ is even} \\ 
p_T/q_T = p^{\prime\prime}/q^{\prime\prime} \, , \text{and,} \, \, 
p_S/q_S = p^{\prime}/q^{\prime} \quad \text{if $\bar{n}$ is odd.}
\end{cases}
\end{equation}
Whereas if $\varepsilon_0 \dots \varepsilon_{\bar{n}-1}=-1$ we have the
symmetric case, namely 
\begin{equation}
\label{eq:resgauss}
\begin{cases}
p_T/q_T = p^{\prime\prime}/q^{\prime\prime} \, , \text{and,} \, \,
p_S/q_S = p^{\prime}/q^{\prime}  \quad \text{if $\bar{n}$ is even} \\ 
p_T/q_T = p^{\prime}/q^{\prime} \, , \text{and,} \, \, 
p_S/q_S = p^{\prime\prime}/q^{\prime\prime} \quad \text{if $\bar{n}$ is odd.}
\end{cases}
\end{equation}
\end{proposition}

\proof
By~\eqref{eq:pnqn2},~\eqref{eq:ptqt} and~\eqref{eq:pmqm} we have:
\begin{equation*}
\begin{cases}
    p_T&=a_{\bar{n}-1}p_{\bar{n}-2}+\varepsilon_{\bar{n}-2}p_{\bar{n}-3}\\
    q_T&=a_{\bar{n}-1}q_{\bar{n}-2}+\varepsilon_{\bar{n}-2}q_{\bar{n}-3} 
  \end{cases} \text{ and }
  \begin{cases}
 p_S &=(a_{\bar{n}}-1)p_{\bar{n}-1}+\varepsilon_{\bar{n}-1}p_{\bar{n}-2}\\
 q_S &=(a_{\bar{n}}-1)q_{\bar{n}-1}+\varepsilon_{\bar{n}-1}q_{\bar{n}-2}
  \end{cases}
\end{equation*}
then:
\begin{equation*}
\frac{p_T+p_S}{q_T+q_S}=
\frac{a_{\bar{n}-1}p_{\bar{n}-2}+\varepsilon_{\bar{n}-2}p_{\bar{n}-3}+
(a_{\bar{n}}-1)p_{\bar{n}-1}+\varepsilon_{\bar{n}-1}p_{\bar{n}-2}}{a_{\bar{n}-1}
q_{\bar{n}-2}+\varepsilon_{\bar{n}-2}q_{\bar{n}-3}+(a_{\bar{n}}-1)q_{\bar{n}-1}+
\varepsilon_{\bar{n}-1}q_{\bar{n}-2}}=  
\frac{p_{\bar{n}}}{q_{\bar{n}}}=\frac{p}{q} \, ,
\end{equation*}
where we used the definition of $p/q$ with its finite continued
fraction of length $\bar{n}$. 
Finally:
\begin{equation*}
\frac{p}{q}-\frac{p_T}{q_T}=\frac{p_{\bar{n}}q_{\bar{n}-1}-p_{\bar{n}-1}
q_{\bar{n}}}{q_{\bar{n}} q_{\bar{n}-1}}= 
\frac{(-1)^{\bar{n}+1}\varepsilon_0 \dots
  \varepsilon_{\bar{n}-1}}{q_{\bar{n}} q_{\bar{n}-1}} \, , 
\end{equation*}
and similarly
\begin{equation*}
\frac{p}{q}-\frac{p_S}{q_S}=\frac{p_{\bar{n}}(q_{\bar{n}}-q_{\bar{n}-1})
-(p_{\bar{n}}-p_{\bar{n}-1})q_{\bar{n}}}{q_{\bar{n}}  q_S}= 
\frac{(-1)^{\bar{n}}\varepsilon_0 \dots \varepsilon_{\bar{n}-1}}{q_{\bar{n}}
  q_S} \, , 
\end{equation*}
from which the proof follows easily.
\endproof

\subsection{The $1/2$--Monoid}
\label{ssec:12monoid}

In this paragraph we will study the monoid $\hat{\mathcal{M}}$ of $GL(2,\Z)$,
 introduced in~\eqref{eq:monoidr} and used in the construction of the $1/2$--Complex
 Bruno function. Our aim is to show its relation with the nearest integer continued 
fraction: for this reason we call it {\em $1/2$--Monoid}. We will prove that given
$p/q \in [0,1)$ we can ``fill'' the matrix 
$g_*=\left(\begin{smallmatrix}p_* &p \\ q_* &q \end{smallmatrix}\right)$ 
in exactly two ways, such that it belongs to $\hat{\mathcal{M}}$ ``following the
nearest integer continued fraction development''.

\begin{proposition}
\label{prop:porp2p6}
Let $p/q \in [0,1)$, $\bar{n}\geq 1$ and assume 
$p/q =[(a_0,\varepsilon_0),\dots,(a_{\bar{n}},\varepsilon_{\bar{n}})]$
 to be the finite continued fraction of $p/q$. 
We claim that the matrices $g_T=
\left(\begin{smallmatrix}p_T &p \\ q_T &q \end{smallmatrix}\right)$ and $g_S=
\left(\begin{smallmatrix}p_S &p \\ q_S &q \end{smallmatrix}\right)$, where the 
 rational $p_T/q_T$ and $p_S/q_S$ have been defined in
 Proposition~\ref{prop:farey2}, 
 are given by:
\begin{align}
\label{eq:12dec}
  g_T&=\hat{\varepsilon}_0g(\hat{a}_1)\dots g(\hat{a}_{\bar{n}-1}) 
\hat{\varepsilon}_{\bar{n}-1}g(a_{\bar{n}})& \quad \text{(Type T)}\\
\label{eq:12decc}
  g_S&=\hat{\varepsilon}_0g(\hat{a}_1)\dots
g(\hat{a}_{\bar{n}-1}) 
\hat{\varepsilon}_{\bar{n}-1}g(a_{\bar{n}}-1)g(1)& \quad \text{(Type S)}\, ,
\end{align}
where for $i=0,\dots,\bar{n}-1$, matrices $\hat{\varepsilon}_i$ and
integer $\hat{a}_i$ are defined by: 
\begin{equation}
\label{eq:epsilon}
(\hat{a}_i,\hat{\varepsilon}_i)=\begin{cases}
(a_i,1) &\text{ if } \varepsilon_i=+1 \\
(a_i-1,\sigma) &\text{ if } \varepsilon_i=-1 \end{cases}\, .
\end{equation}
\end{proposition}

Before to prove the Proposition we make the following:
\begin{remark}
\label{rem:2tipes}
It results $\hat{a}_i\geq 2$ for all $i$, in fact whenever
$\varepsilon_i=-1$  
one has $a_i \geq 3$ (see Remark~\ref{rem:conv}). Because $p/q \in
[0,1)$ the first couple $(a_0,\varepsilon_0)$ can 
only be one of the following two: $(0,+1)$ if $p/q\in [0,1/2]$ or $(1,-1)$ if 
$p/q\in (1/2,1)$.
\end{remark}
\proof
Let $k \leq \bar{n}$ and let us introduce matrices 
$\hat{\varepsilon}_0,\dots,\hat{\varepsilon}_{k-1}$ and integers
$\hat{a}_0,\dots,\hat{a}_{k}$ as in~\eqref{eq:epsilon} according to the 
continued fraction of $p/q$. Then we claim that
$g=\hat{\varepsilon}_0g(\hat{a}_1)\dots
g(\hat{a}_{k-1})\hat{\varepsilon}_{k-1}g(a_k)$ is equal to
$\left(\begin{smallmatrix}p_{k-1} &p_{k} \\ q_{k-1} &q_k
\end{smallmatrix}\right)$ 
where $p_k/q_k=[(a_0,\varepsilon_0),\dots,(a_k,\varepsilon_k)]$. This
can be proved by induction (use Remark~\ref{rem:2tipes} to prove the
basis of induction) and then~\eqref{eq:12dec} follows by putting $k=\bar{n}$. 
To prove~\eqref{eq:12decc} it is enough to calculate:
\begin{equation*}
\hat{\varepsilon}_0g(\hat{a}_1)\dots
g(\hat{a}_{k-1})\hat{\varepsilon}_{k-1}g(a_k-1)g(1)= 
\left(\begin{matrix}p_{k-1} &p_{k}-p_{k-1} \\ q_{k-1} &q_{k}-q_{k-1}
\end{matrix}\right) 
\left(\begin{matrix}0 &1 \\ 1 &1 \end{matrix}\right)=
\left(\begin{matrix}p_{k}-p_{k-1} &p_k \\ q_{k}-q_{k-1}
  &q_k\end{matrix}\right) \, . 
\end{equation*}
\endproof

\begin{remark}
Clearly matrices of type T belong to $\Hat{\mathcal{M}}$ (because $a_{\bar{n}}\geq 2$)
, whereas those of type S belong to the monoid if and only if the continued fraction of the
rational $p/q$ ends with a couple $(a_{\bar{n}},\varepsilon_{\bar{n}})=(2,1)$,
in fact in this way the matrix $g_S$ is given by
$\hat{\varepsilon}_0g(\hat{a}_1)\dots g(\hat{a}_{\bar{n}-1}) 
g(1)g(1)=\hat{\varepsilon}_0g(\hat{a}_1)\dots g(\hat{a}_{\bar{n}-1})
\sigma g(2)$, where we used that $\varepsilon_{\bar{n}-1}=1$ (because
$a_{\bar{n}}=2$)  
and $\sigma g(m)=g(m-1)g(1)$ for all $m\geq 2$.

Remark also that, if $g$ is of type T then it cannot ends with $\sigma g(2)$, 
in fact this will implies a continued fraction ending with  
$[\dots,(a_{\bar{n}-1},-1),(2,1)]$, but we know that this is impossible and so:
or $\varepsilon_{\bar{n}-1}=-1$, $a_{\bar{n}}\geq 3$ and $g_T=\ldots
\sigma g(a_{\bar{n}})$, 
either $\varepsilon_{\bar{n}-1}=+1$, $a_{\bar{n}}\geq 2$  and $g_T=\ldots
g(a_{\bar{n}-1})g(a_{\bar{n}})$. 
\end{remark}
With the following proposition we will prove that $\Hat{\mathcal{M}}$
is the union 
of matrices of type T and of type S with
$(a_{\bar{n}},\varepsilon_{\bar{n}})=(2,1)$. 
Let us denote by $\mathcal{M}_T$ the monoid of matrices of type T and
$\mathcal{M}_S$ those 
of type S, with $(a_{\bar{n}},\varepsilon_{\bar{n}})=(2,1)$.

\begin{proposition}[The $1/2$--Monoid]
\label{prop:12monoidTS}
$\Hat{\mathcal{M}}=\mathcal{M}_T\cup \mathcal{M}_S$.
\end{proposition}
\proof
Clearly $\mathcal{M}_T\cup \mathcal{M}_S \subset
\Hat{\mathcal{M}}$. Let us prove 
the other inclusion. Let $r\geq 1$, $m_1,\dots,m_r \geq 2$, 
$\hat{\varepsilon}_0 ,\dots,\hat{\varepsilon}_{r-1} \in \{ 1, \sigma
\}$, such that 
$g=\hat{\varepsilon}_0 g(m_1)\dots\hat{\varepsilon}_{r-1}
g(m_r)\in\Hat{\mathcal{M}}$. 

Let us consider two cases, first one: $\hat{\varepsilon}_{r-1}=\sigma$ and
$m_r \geq 3$ or $\hat{\varepsilon}_{r-1}=1$ and $m_r \geq 2$; second
case: $\hat{\varepsilon}_{r-1}=\sigma$ and
$m_r =2$. In the former case to $g$ we associate a continued fraction
by introducing, for $i=1,\dots,r-1$:
\begin{equation*}
(a_0,\varepsilon_0)=\begin{cases}
(0,+1) &\text{ if }\hat{\varepsilon}_i=1 \\
(1,-1) &\text{ if }\hat{\varepsilon}_i=\sigma
\end{cases}
\, ,\quad
(a_i,\varepsilon_i)=\begin{cases}
(m_i,+1) &\text{ if }\hat{\varepsilon}_i=1 \\
(m_i+1,-1) &\text{ if }\hat{\varepsilon}_i=\sigma
\end{cases}
\, , \quad a_r=m_r \, .
\end{equation*}
$[(a_0,\varepsilon_0),\dots,(a_r,\varepsilon_r)]$
represents some rational $p/q$, let us define as before $p_T/q_T$ and 
then $g_T=\left(\begin{smallmatrix}p_T & p \\ q_T & q
\end{smallmatrix}\right)=
\hat{\varepsilon}_0g(\hat{a}_1)\dots\hat{\varepsilon}_{r-1}g(\hat{a}_r)$
(by Proposition~\ref{prop:porp2p6}, where we also defined $\hat{a}_i$'s). 
Observe that $\hat{a}_i=m_i$ to conclude $g=g_T\in\mathcal{M}_T$.

The second case can be treated similarly. Now to $g$ we associate the
continued fraction $[(a_0,\varepsilon_0),\dots,(a_{r-1},1),(2,1)]$ 
where, for $i=1,\dots,r-2$:
\begin{equation*}
(a_0,\varepsilon_0)=\begin{cases}
(0,+1) &\text{ if }\hat{\varepsilon}_i=1 \\
(1,-1) &\text{ if }\hat{\varepsilon}_i=\sigma
\end{cases}
\, ,\quad
(a_i,\varepsilon_i)=\begin{cases}
(m_i,+1) &\text{ if }\hat{\varepsilon}_i=1 \\
(m_i+1,-1) &\text{ if }\hat{\varepsilon}_i=\sigma
\end{cases}
\, , \, a_{r-1}=m_{r-1} \, .
\end{equation*}
Let $[(a_0,\varepsilon_0),\dots,(a_{r-1},1),(2,1)]$ be some
rational $p/q$, define as before 
$p_S/q_S=[(a_0,\varepsilon_0),\dots,(a_{r-1},1),(1,1)]$, 
then by Proposition~\ref{prop:porp2p6}:
\begin{equation*}
g_S=\hat{\varepsilon}_0g(\hat{a}_1)\dots g(\hat{a}_{r-1})g(1)g(1)=
\hat{\varepsilon}_0g(\hat{a}_1)\dots g(\hat{a}_{r-1})\sigma g(2)=g
\end{equation*}
and it belongs to $\mathcal{M}_S$.
\endproof

To end this section we introduce a third characterization of the
$1/2$--Monoid, which corrects a small error in~\S A.4.4 page 836 of~\cite{MMYc},
and which will be useful to construct the numerical algorithm for the
$1/2$--complex Bruno function.
\begin{proposition}
Let $g=\left(\begin{smallmatrix} a & b \\ c &d\end{smallmatrix}\right)\in G$.
Then $g$ belongs to $\Hat{\mathcal{M}}$ if and only if
$d\geq b>0$, $c\geq a\geq 0$ and $d\geq \mathcal{G} c$, where 
$\mathcal{G}=(\sqrt{5}+1)/2$.
\end{proposition}
The proof can be done by direct computation and we omit it. We end this 
part with the following:

\begin{remark}[The Gauss Monoid]
\label{rem:gaussmonoid}
In~\cite{MMYc} authors considered the complex Bruno function constructed
using the Monoid $\mathcal{M}$:
\begin{equation*}
\mathcal{M}=\Big \{ g=\left(\begin{smallmatrix} a & b \\ c &d\end{smallmatrix}\right)\in G:
d\geq b\geq a \geq 0\text{ and } d\geq c\geq a \Big \} \, .
\end{equation*}
We recall that according to the Gauss continued fraction algorithm we always have $\varepsilon_l=+1$,
we can then prove modified version of Propositions~\ref{prop:farey2} and~\ref{prop:porp2p6}
to conclude that $\mathcal{M}$ is constructed ``following'' the Gauss continued
fraction algorithm: starting from $p/q \in (0,1)$, we complete the matrix
$g_*=\left(\begin{smallmatrix} p_* & p \\ q_* & q \end{smallmatrix}\right)$ into
$g_S$ and $g_T$, where $p_S/q_S$ and $p_T/q_T$ are obtained with the 
Truncate and Subtract operations acting on the Gauss finite continued fraction of $p/q$.
\end{remark}

\subsection{An Algorithm for the $1/2$--complex Bruno function}
\label{ssec:cmplxbruno}

Using the results of the previous sections we are now able to give an algorithm to
compute the $1/2$--Complex Bruno function. Let us rewrite definition~\eqref{eq:defcomplxB}:
\begin{equation*}
\B(z)=\sum_{n\in \Z}\left[\sum_{g\in\hat{\mathcal{M}}}L_g\left(1+L_{\sigma}\right)\right] 
\varphi_{1/2}(z-n) \, , 
\end{equation*}
where $\varphi_{1/2}(z)=-\frac{1}{\pi}\li{\frac{1}{2z}}+
\frac{1}{\pi}\log 2 \log \left(1-\frac{1}{2z} \right)$ and the action $L_g$ has been 
defined in~\eqref{eq:Laction}. From the previous sections we know that the sum over 
$\hat{\mathcal{M}}$ can be replaced by a sum over $p/q \in [0,1)$ , $(p,q)=1$, in such
a way that to each $p/q$ we associate the matrix $g_T$, and also $g_S$ whenever the continued
fraction of $p/q$ ends with $(a_{\bar{n}},\varepsilon_{\bar{n}})=(2,+1)$.

Using the periodicity and the parity properties of $\B$ we can restrict to $\Re z \in [0,1/2]$. 
Let us consider the contribution of some $p/q \in [0,1)$ to $\B$. Because of the form of $\varphi_{1/2}$ and of the action $L_g$ we remark that the
larger is the denominator of the fraction, the smaller is its contribution to the sum; moreover different rational numbers with the same denominator give
comparable contributions, so we decide to order the rationals w.r.t. increasing
denominators: in other words {\em according to the Farey Series}. 
A similar statement holds w.r.t. the sum over $\Z$: large $n$'s 
give small contribution to the sum. We then introduce two {\em cut--off}
 to effectively compute~\eqref{eq:defcomplxB}: $N_{max}$ denoting the largest order
of the Farey Series considered and $k_1$ the largest (in modulus) $n\in\Z$ which 
contributes to the sum over integers~\footnote{
For technical reasons, we prefer to introduce a third cut--off, $k_2$. We refer
to Appendix~\ref{sec:numcons} to explain the role of this cut--off.}.

\begin{center}
  \begin{figure}[hb]
    \begin{center}
    \mbox{\subfigure
    {\includegraphics[scale=0.3,angle=-90]{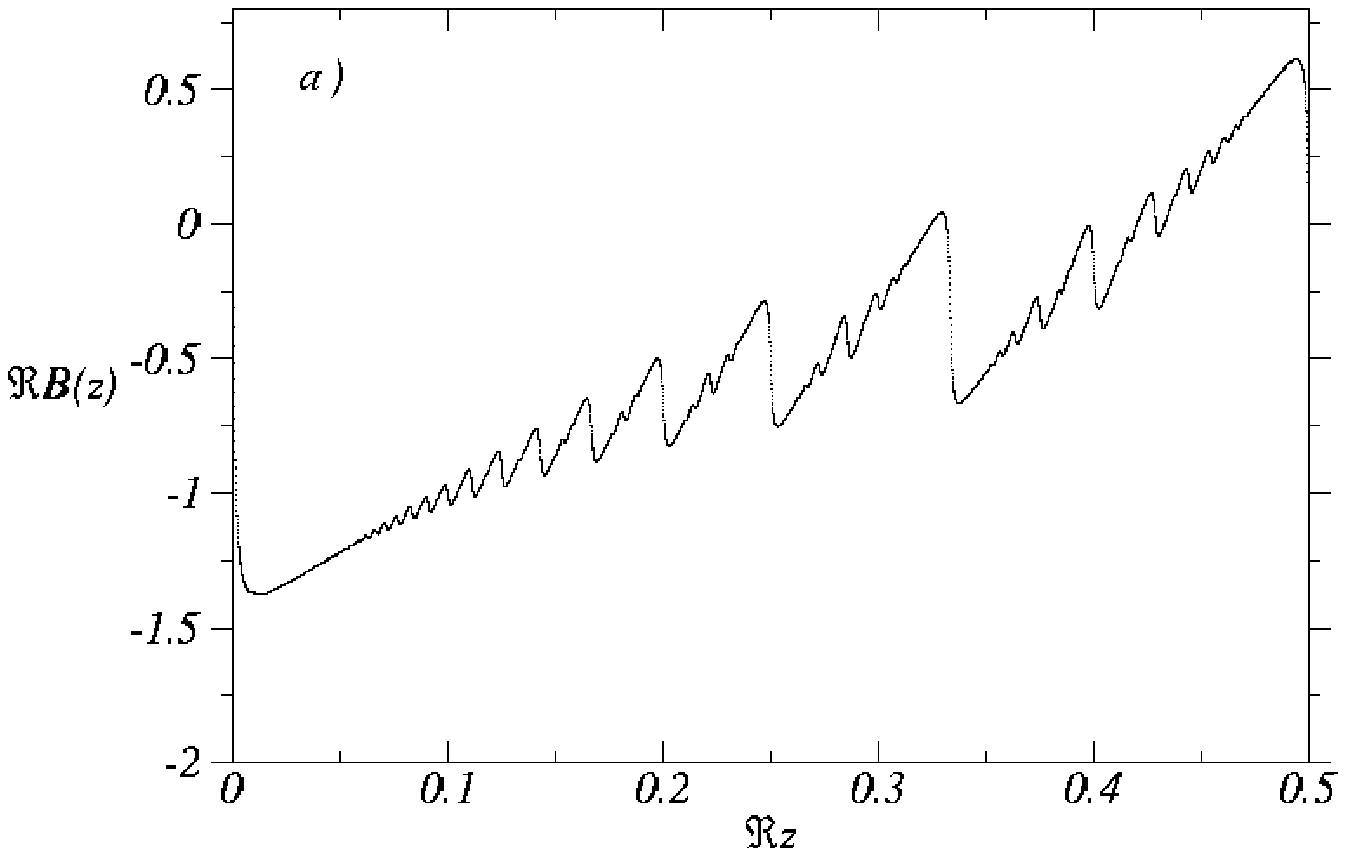}}\quad\quad
    \subfigure
    {\includegraphics[scale=0.3,angle=-90]{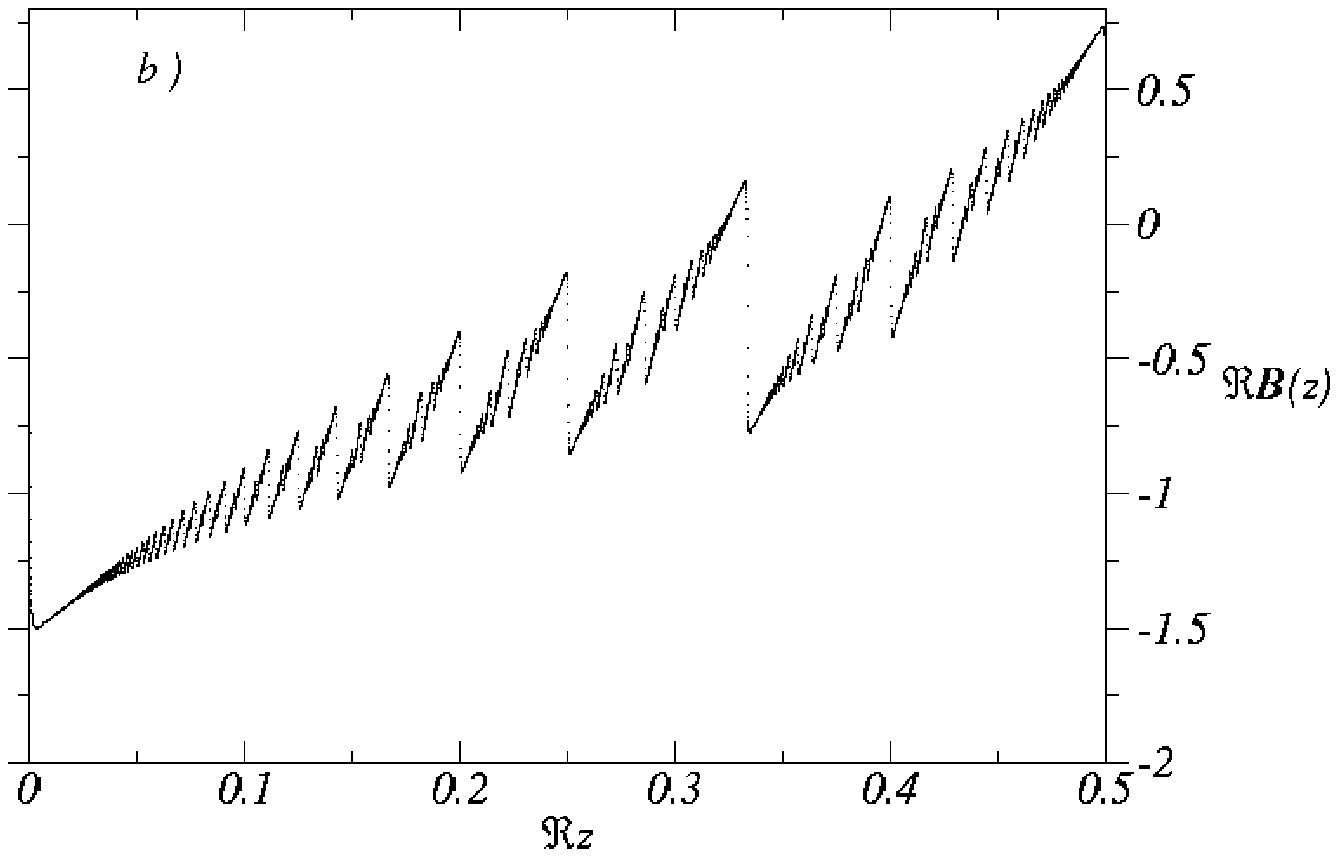}}
       }
    \mbox{\subfigure
    {\includegraphics[scale=0.3,angle=-90]{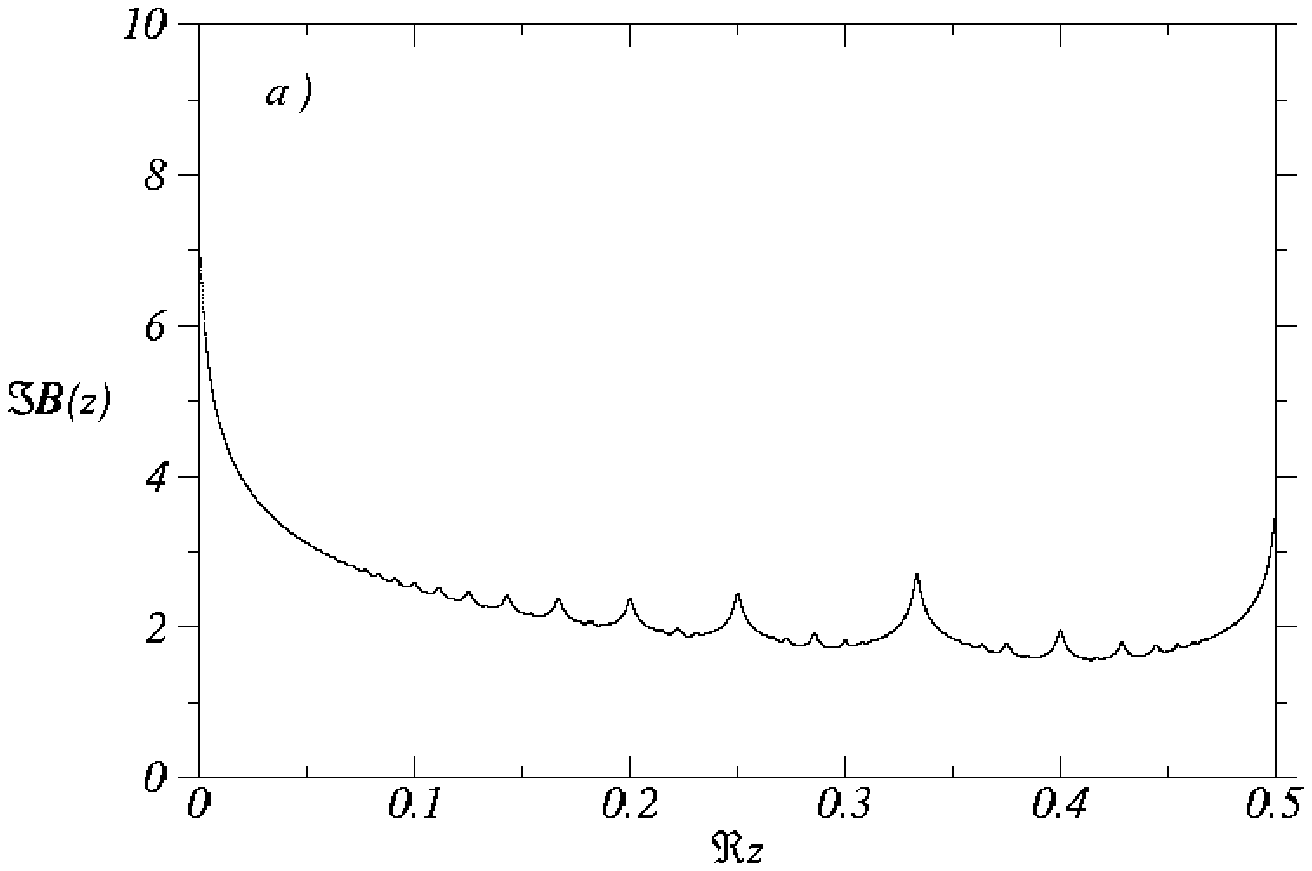}}\quad\quad
    \subfigure
    {\includegraphics[scale=0.3,angle=-90]{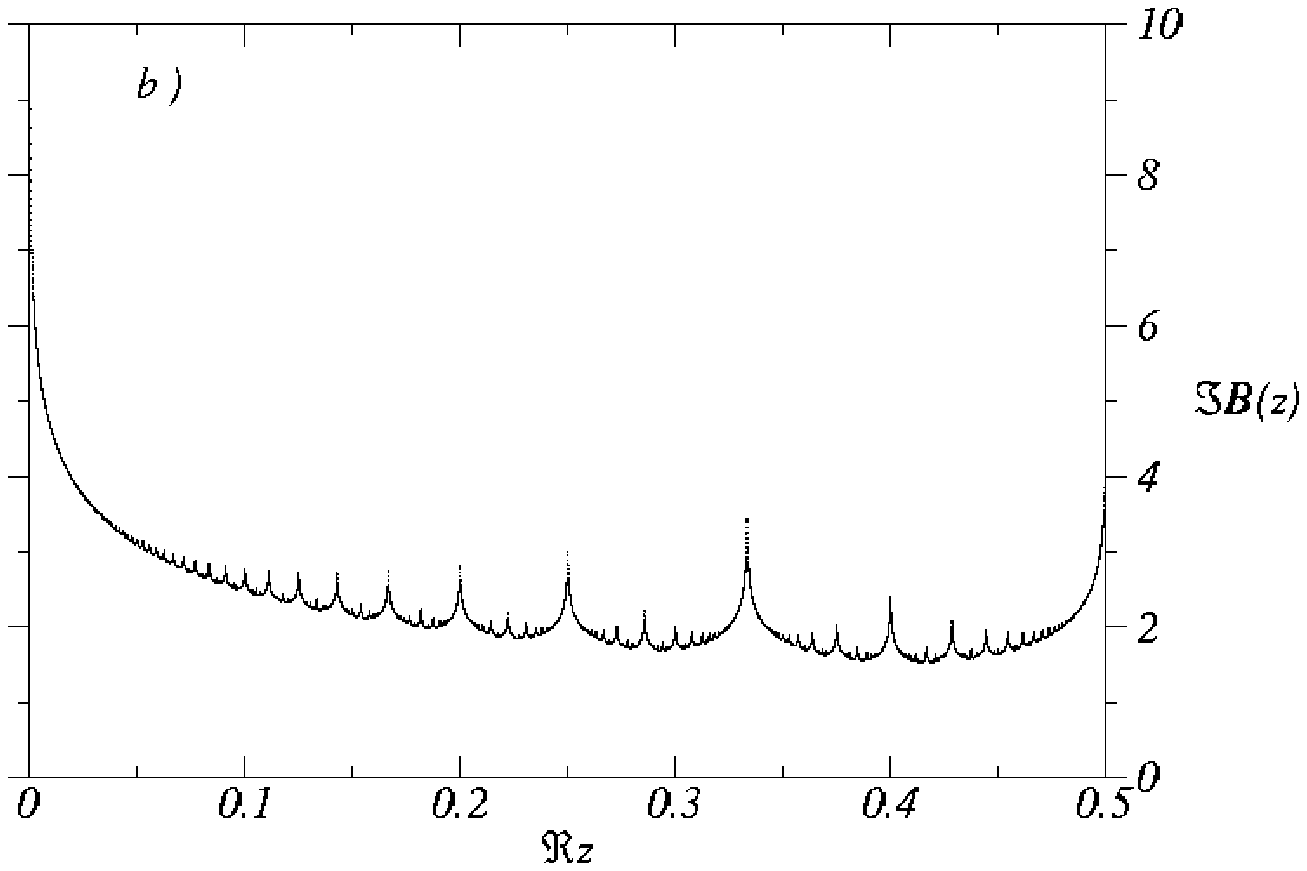}}
       }
    \end{center}
   \caption{Plot of $\B(z)$ vs $\Re z$ at $\Im z$ fixed. 
    The top line contains $\Re \B$ whereas on the bottom line we plot
    $\Im \B$. The column a) is for $\Im z= 10^{-3}$ whereas column b)
    is for $\Im z= 10^{-4}$. 
     Each plot has $10000$ points $\Re z$ 
   uniformly distributed in $[0,1/2]$. $k_1=80$, $k_2=20$, $N_{max}=151$.} 
  \label{fig:imbrunoc}
  \end{figure}
\end{center}

\begin{center}
  \begin{figure}[ht]
    \begin{center}
    \mbox{\subfigure
    {\includegraphics[scale=0.3,angle=-90]{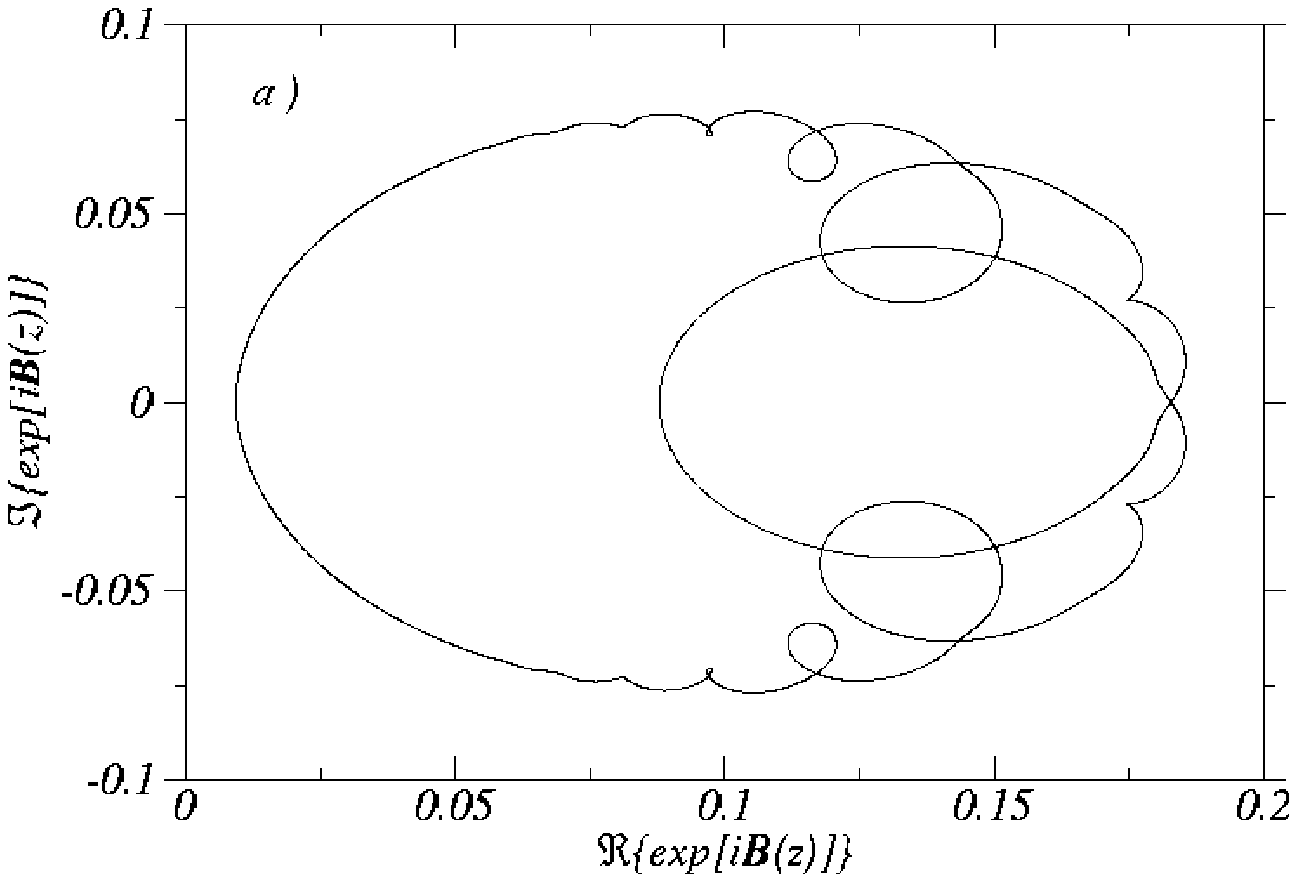}}\quad\quad
    \subfigure
    {\includegraphics[scale=0.3,angle=-90]{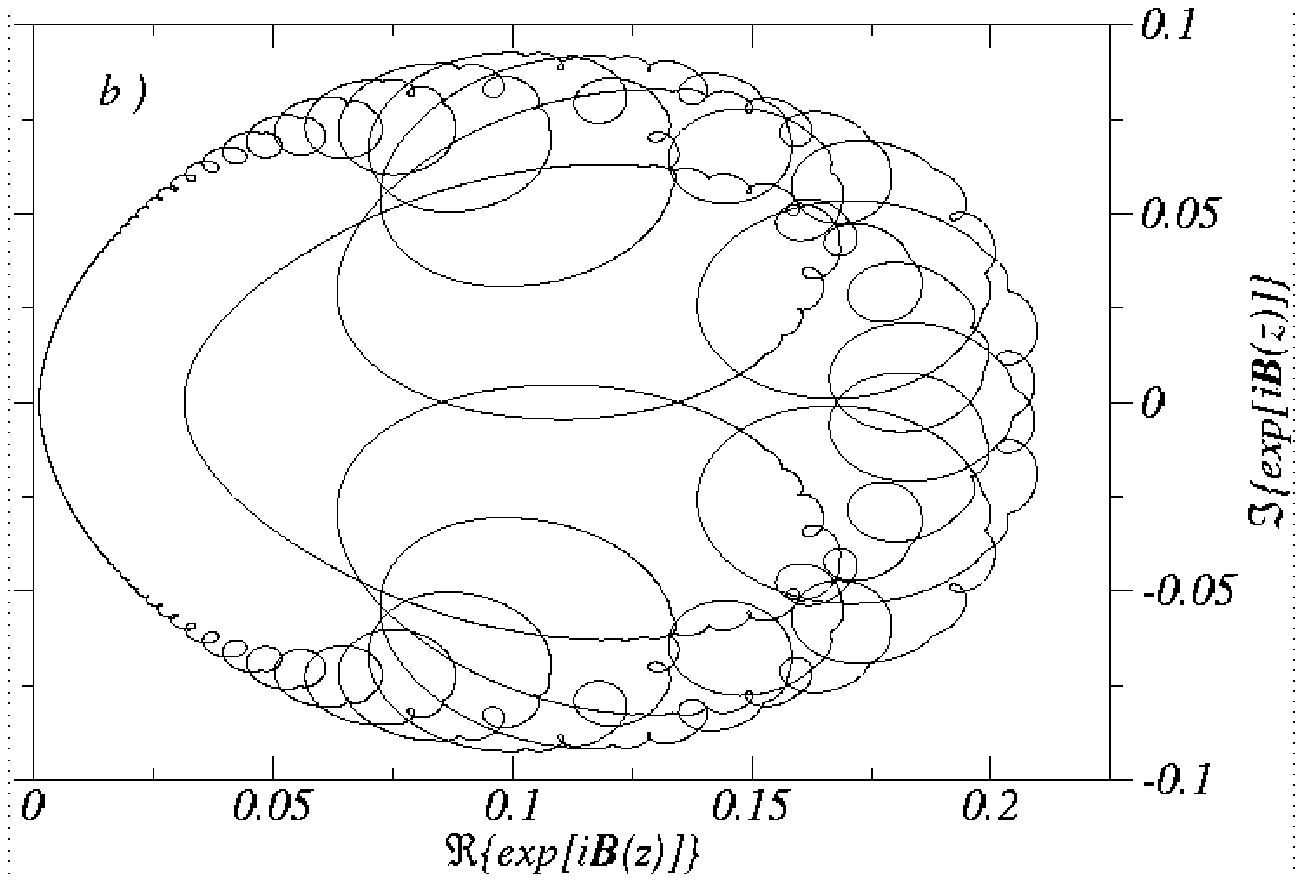}}
    }
    \end{center}
  \caption{Polar plot of $e^{i\B\left( z \right)}$ for fixed values of $\Im z$.
  a)~$\Im z= 10^{-2}$ and b)~$\Im z= 10^{-3}$. $12000$ points
  uniformly distributed in $[0,1]$, $k_1=80$, $k_2=20$, $N_{max}=151$.} 
  \label{fig:polareiB}
  \end{figure}
\end{center}

Then the $1/2$--complex Bruno function can be numerically approximated by:
\begin{equation}
\label{eq:numapprx}
\B(z) \sim \sum_{|n|\leq k_1}\sum^{\prime}_{p/q \in \mathcal{F}_{N_{max}}}
L_{\left( \begin{smallmatrix} p_* & p \\ q_* & q\end{smallmatrix}\right)}
\left(1+L_{\sigma}\right)\varphi_{1/2}(z-n) \, , 
\end{equation}
where $p_*/q_*\in \{ p_T/q_T,p_S/q_S \}$ and the sum is restricted to fractions
such that $q \geq \mathcal{G}q_S$, (being always $q \geq \mathcal{G}q_T$).
 This approximation can be made as precise
 as we want, by choosing $N_{max}$ and $k_1$ large enough, in fact~\eqref{eq:defcomplxB}
 can be obtained as double limit $N_{max}\rightarrow +\infty$ and
 $k_{1}\rightarrow +\infty$. In Appendix~\ref{sec:numcons} we will give
numerical results showing the convergence of~\eqref{eq:numapprx} varying the cut-off
values, the convergence of $\Im \B(z)$ to $B(\Re z)$ when $\Im z \rightarrow 0$ and
$\Re z \in \mathcal{B}$, and the $\pi/q$--jumps of $\Re \B(z)$ when $z\rightarrow p/q$,
 as proved in~\cite{MMYc}. In Figure~\ref{fig:imbrunoc} we show some
 plots of $\B(z)$ for fixed (small) values of
$\Im z$ and $\Re z \in [0,1/2]$, whereas in Figure~\ref{fig:polareiB} we show two polar plots
of $e^{i\B\left(z\right)}$.

\section{The Yoccoz Function}
\label{sec:yoccozf}

The aim of this section is to briefly introduce the algorithm used to compute the Yoccoz
Function, $U\!(\lambda)$, introduced in \S~\ref{ssec:yoccozfunc}. Let $\lambda \in \D^*$, 
let $P_{\lambda}(z)=\lambda z \left( 1-z\right)$ be the quadratic polynomial and let
us introduce the polynomials: $U_n(\lambda)=\lambda^{-n}
P^{\circ n}_{\lambda}\left(1/2\right)$. Then we recall that the Yoccoz function 
is the uniform limit, over compact subsets of $\D$, of $U_n(\lambda)$.

From~\eqref{eq:siegel} and its original definition, 
$H_{\lambda}\left(U\!\left(\lambda\right)\right)=1/2$, we get: 
\begin{equation}
\label{eq:un1}
U\!\left(\lambda\right)=
\lambda^{-n}H^{-1}_{\lambda}\left(\lambda^n
U_n\!\left(\lambda\right)\right)\, , 
\end{equation}
for all integer $n$. Hence to compute $U\!\left(\lambda\right)$ we need to know
how close is $H^{-1}_{\lambda}$ to the identity, near zero and this
can be done using 
 some standard distortion estimates~\cite{BHH}. So for any fixed
 $\lambda \in \D^*$, 
we can find $n=n\left(\lambda\right)$ s.t. $P^{\circ
  n}_{\lambda}\left(1/2\right)$  
is contained in some fixed disk on which we can apply the distortion
estimate and then 
from~\eqref{eq:un1} compute an approximation to $U\!\left(\lambda\right)$
with a prescribed precision $\epsilon_{U}$.

\begin{remark}[Parity of Yoccoz's function]
\label{rem:parityfourcoeff}
Let us observe the following facts. Assume $\lambda = e^{2\pi i
  \left(x+it\right)}$, with $t>0$ fixed, and $x$ varying in $(0,1/2)$ and let us
  introduce $u(x)=U\!\left(e^{2\pi i \left(x+it\right)}\right)$, to stress the
  dependence on $x$ only. Then we claim that $\Re u(-x)=\Re u(x)$ and
  $\Im u(-x)=-\Im u(x)$. The proof can be done as follows. First remark
  that $\lambda$, as a function of $x$, is mapped into
  $\bar{\lambda}$, when $x\mapsto -x$; then is enough to observe that
  polynomials $U_n\left(\lambda \right)$ verify, for $n\geq 2$,
  $U_n\left( \bar{\lambda}\right)=\overline{U_n\left(\lambda\right)}$,
  namely:
  \begin{equation*}
    \Re U_n\!\left(e^{2\pi i \left(x+it\right)}\right)=\Re 
    U_n\!\left(e^{2\pi i \left(-x+it\right)}\right)
\text{ and }  \Im U_n\!\left(e^{2\pi i \left(x+it\right)}\right)=-\Im 
    U_n\!\left(e^{2\pi i \left(-x+it\right)}\right) \, .
  \end{equation*}
A similar statement holds for $\log U \!\left( \lambda \right)$.

Using the $\Z$--periodicity, we consider the Fourier series of 
$U\!\left(e^{2\pi iz}\right)$ and using an argument similar to the one 
of Remark~\ref{rem:fourcoeffB} we conclude that all the Fourier coefficients
are real and zero for negative Fourier modes. Clearly Taylor's coefficients
of $U\!\left(\lambda\right)$ coincide with Fourier coefficients 
of $U\!\left(e^{2\pi iz}\right)$.
\end{remark}

\begin{center}
  \begin{figure}[ht]
    \begin{center}
    \mbox{\subfigure
    {\includegraphics[scale=0.3,angle=-90]{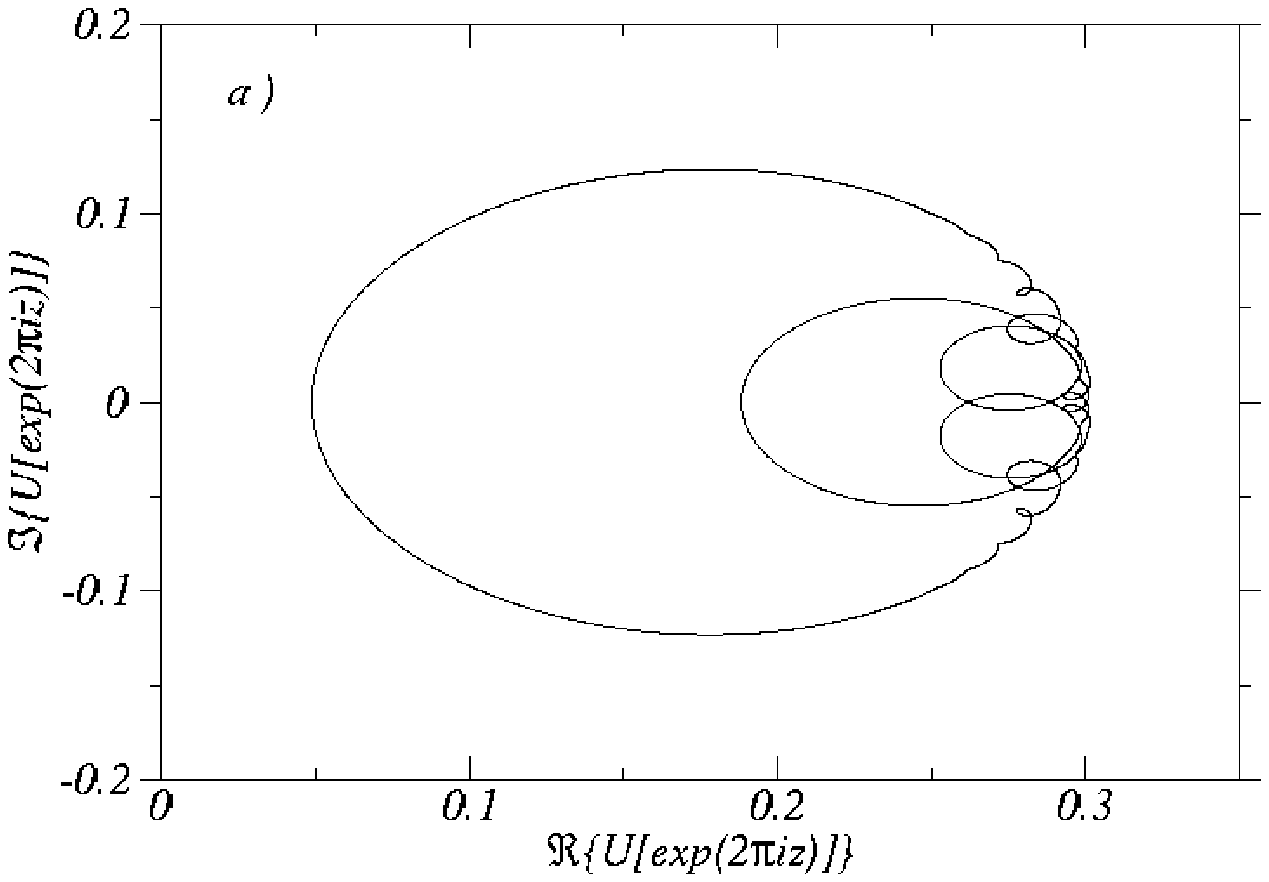}}\quad\quad
    \subfigure
    {\includegraphics[scale=0.3,angle=-90]{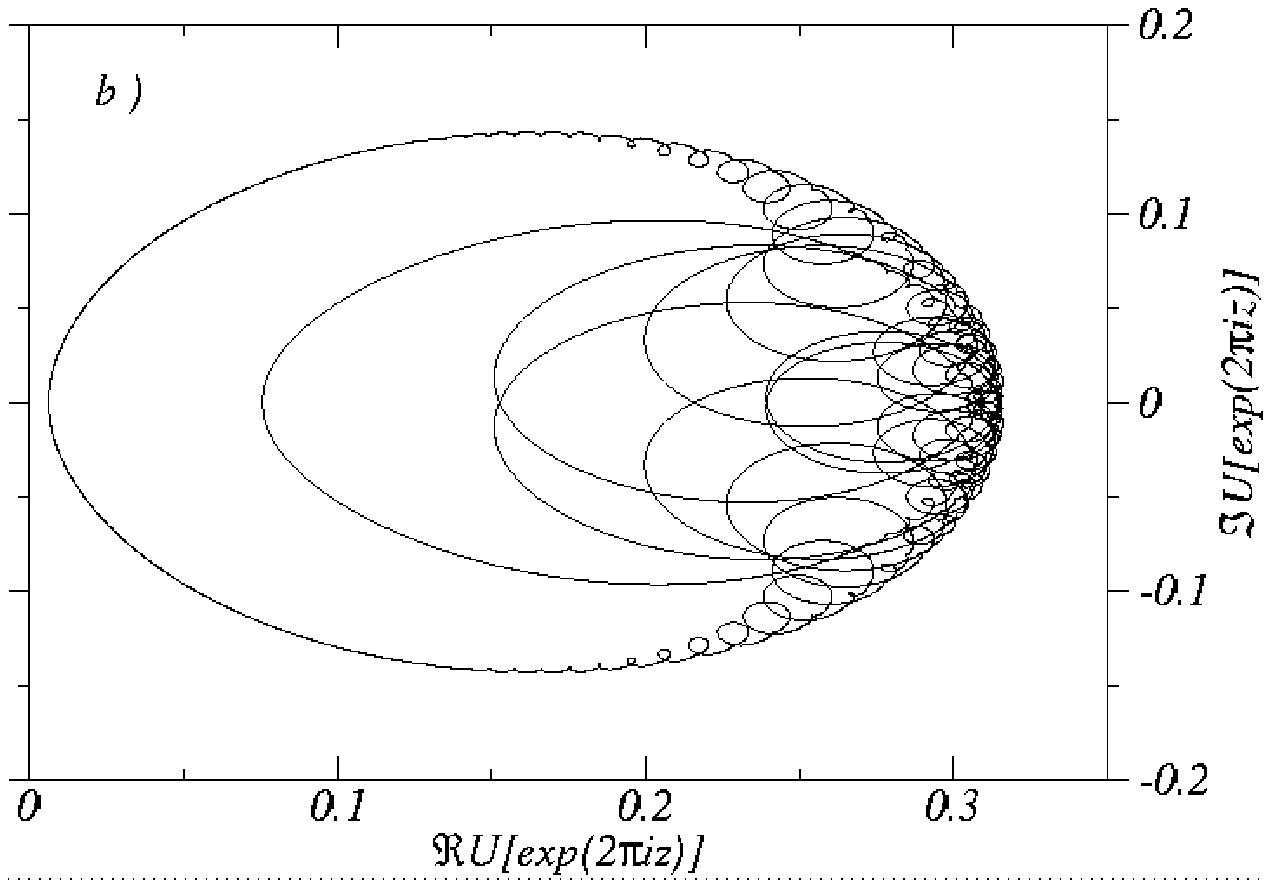}}
    }
    \end{center}
  \caption{Polar plot of $U\!\left( e^{2\pi i z}\right)$ for fixed values of $\Im z$.
  a)~$\Im z= 10^{-2}$ and b)~$\Im z= 10^{-3}$. $12000$ points
  uniformly distributed in $[0,1]$, $\epsilon_{U}=10^{-3}$.} 
  \label{fig:polarU}
  \end{figure}
\end{center}

\begin{center}
  \begin{figure}[ht]
    \begin{center}
    \mbox{\subfigure
    {\includegraphics[scale=0.3,angle=-90]{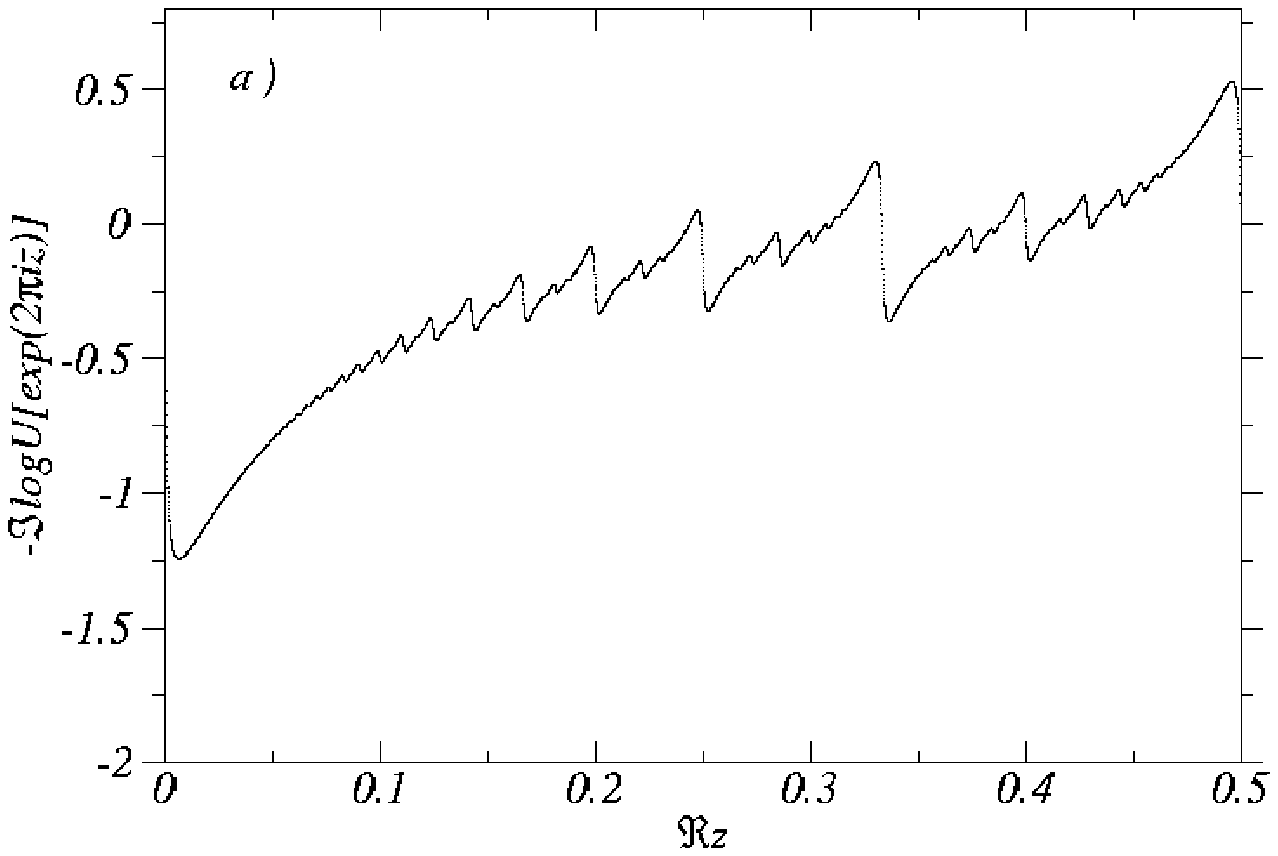}}\quad\quad
    \subfigure
    {\includegraphics[scale=0.3,angle=-90]{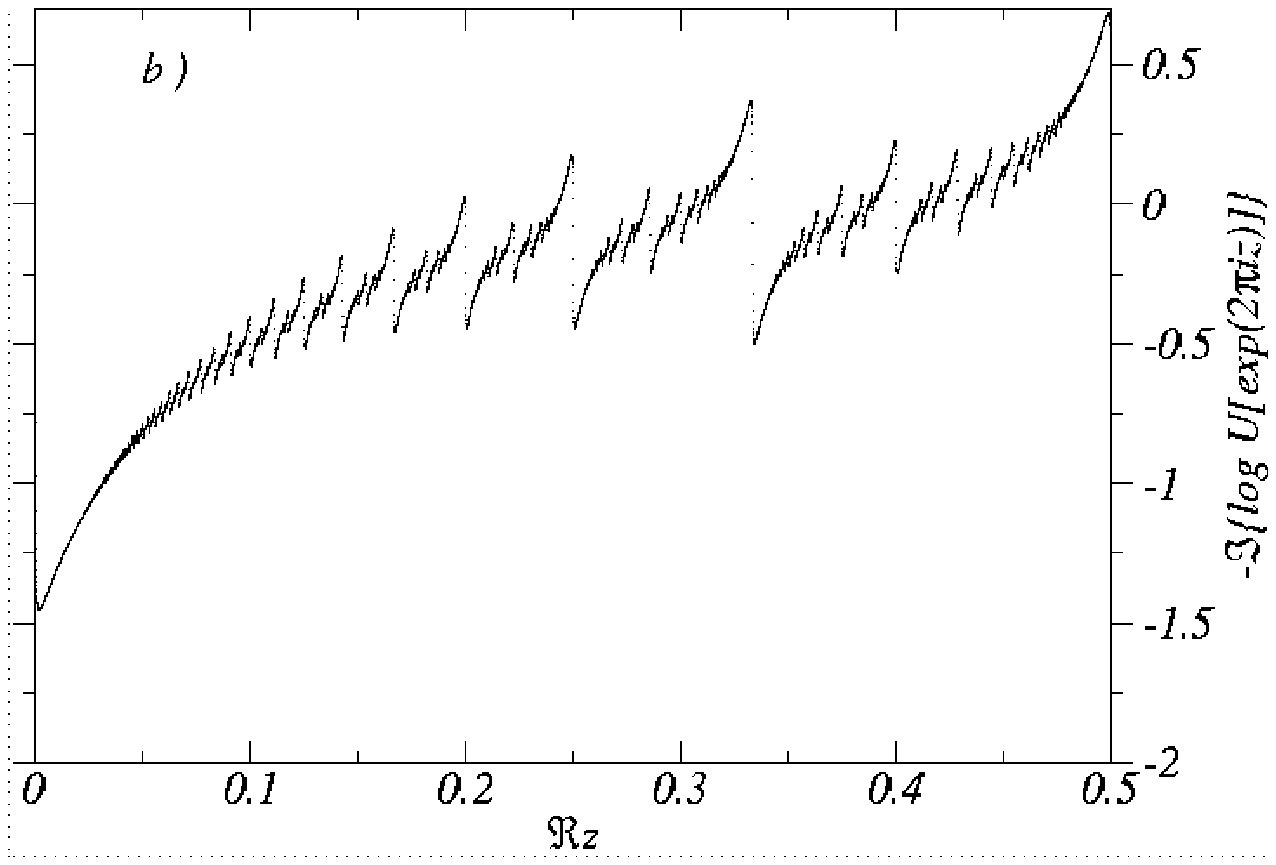}}\quad\quad
       }
    \mbox{\subfigure
    {\includegraphics[scale=0.3,angle=-90]{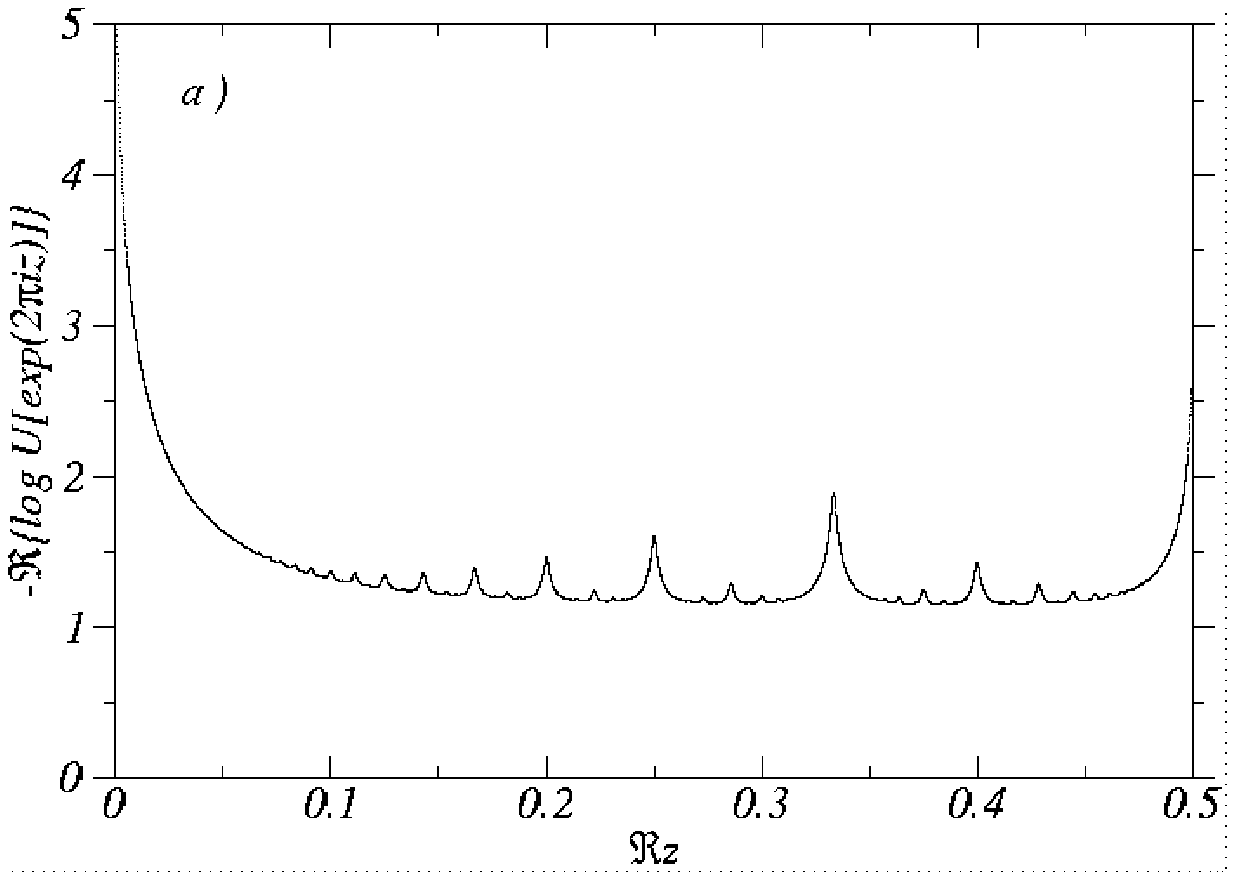}}\quad\quad
    \subfigure
    {\includegraphics[scale=0.3,angle=-90]{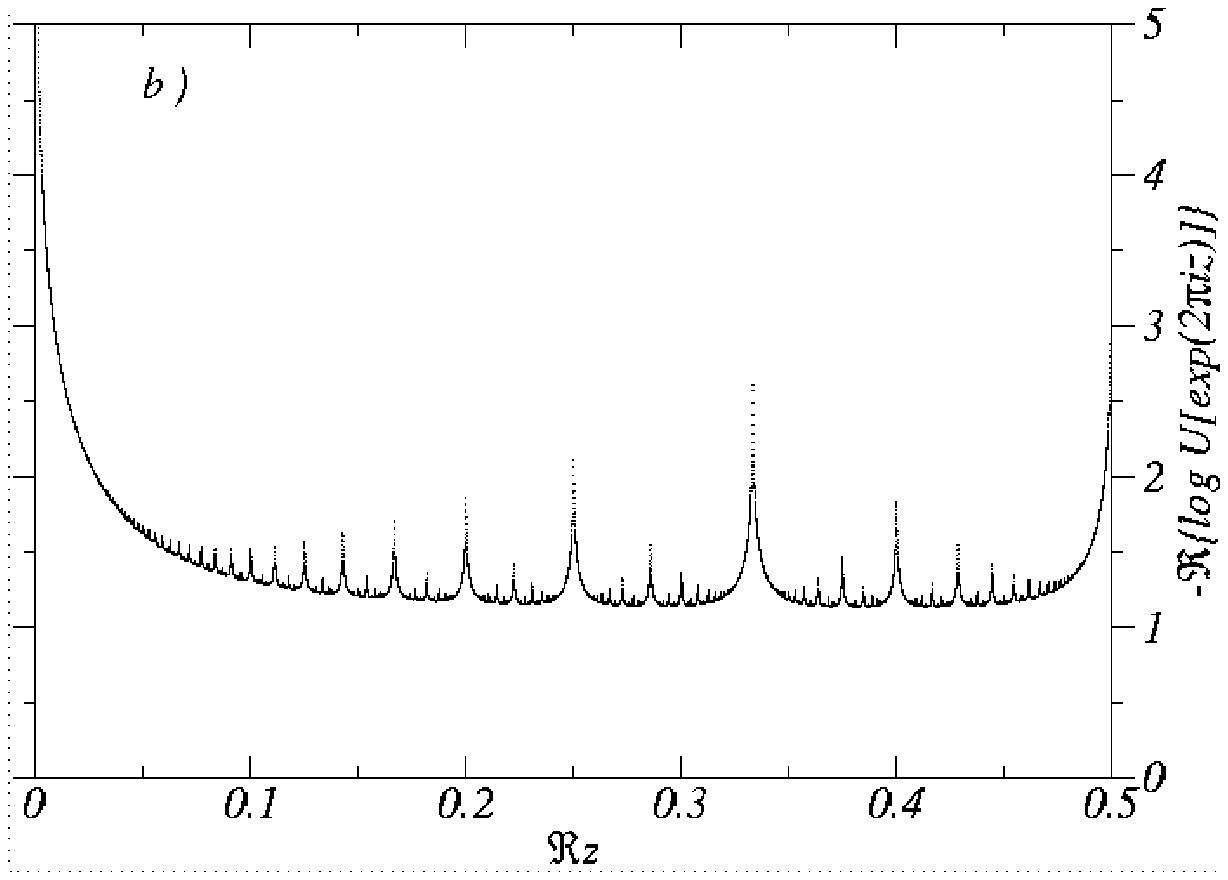}}\quad\quad
       }
    \end{center}
  \caption{Plot of $-\log U\!\left( e^{2\pi i z}\right)$ vs $\Re z$ at fixed
  $\Im z$. On the top we plot the imaginary part whereas on
  the bottom the real part. Column  a) is for~$\Im z= 10^{-3}$, whereas
  on the column b) we show $\Im z= 10^{-4}$
  Each plot has $10000$ points uniformly distributed in
  $[0,1/2]$, $\epsilon_{U}=10^{-3}$.} 
  \label{fig:imreU}
  \end{figure}
\end{center}

Figure~\ref{fig:polarU} shows some polar plots of $U\!\left(
e^{2\pi i z}\right)$, for  
different values of $\Im z >0$, whereas in Figure~\ref{fig:imreU} real and 
imaginary parts of $-\log U\!\left( e^{2\pi i z}\right)$ are
given. Compare with Figures~\ref{fig:polareiB} and
~\ref{fig:imbrunoc}.


Let us conclude this section with the following remark.
\begin{remark}
\label{rem:polinomi}
In Figure~\ref{fig:polarV} we show some polar plots of the ``Yoccoz
function'' used in~\cite{BHH} (Figure 2, page 484):
they don't look like our previous pictures. Here is the reason. 
They take the following
quadratic polynomial $Q_{\lambda}(z)=\lambda z+z^2$, which can be conjugate to our
choice, $P_{\lambda}(z)=\lambda z\left( 1-z\right)$, using $\Lambda(z)=-\lambda z$:
\begin{equation*}
\Lambda \circ P_{\lambda}=Q_{\lambda}\circ \Lambda \, .
\end{equation*}
Let us denote by $V\!\left(\lambda\right)$ the Yoccoz function for the polynomial
 $Q_{\lambda}$, then we claim that:
\begin{equation*}
-\lambda U\!\left(\lambda\right) = V\!\left(\lambda\right) \, ,
\end{equation*}
which explain completely the relation between Figure~\ref{fig:polarU} and 
Figure~\ref{fig:polarV}. Because $-\Im \log U\!\left(\lambda\right)$ exhibits the same jumps at
 rationals as the real part of the complex Bruno function does, we choose
 the quadratic 
polynomial in the form $P_{\lambda}$.
\begin{center}
  \begin{figure}[ht]
    \begin{center}
    \mbox{\subfigure
    {\includegraphics[scale=0.3,angle=-90]{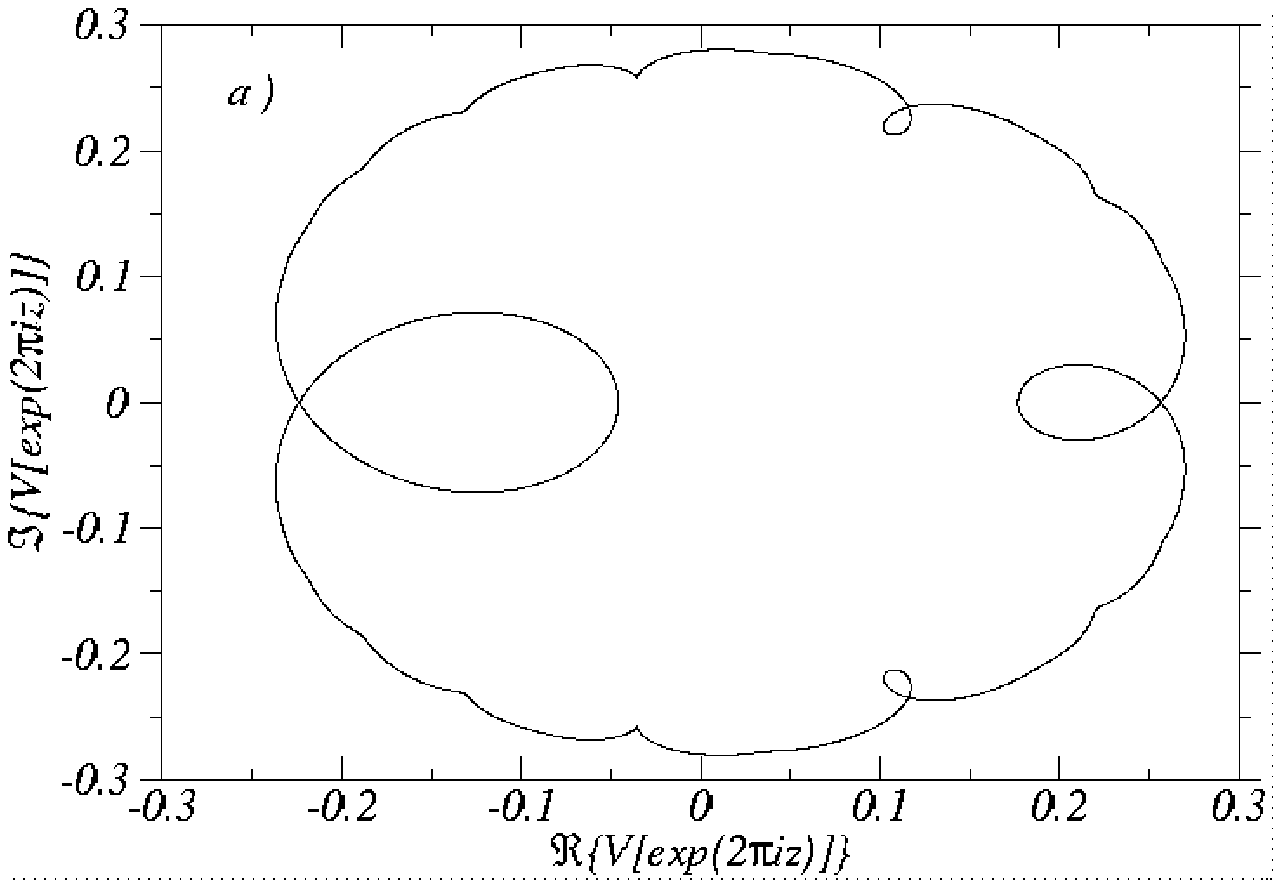}}\quad\quad
    \subfigure
    {\includegraphics[scale=0.3,angle=-90]{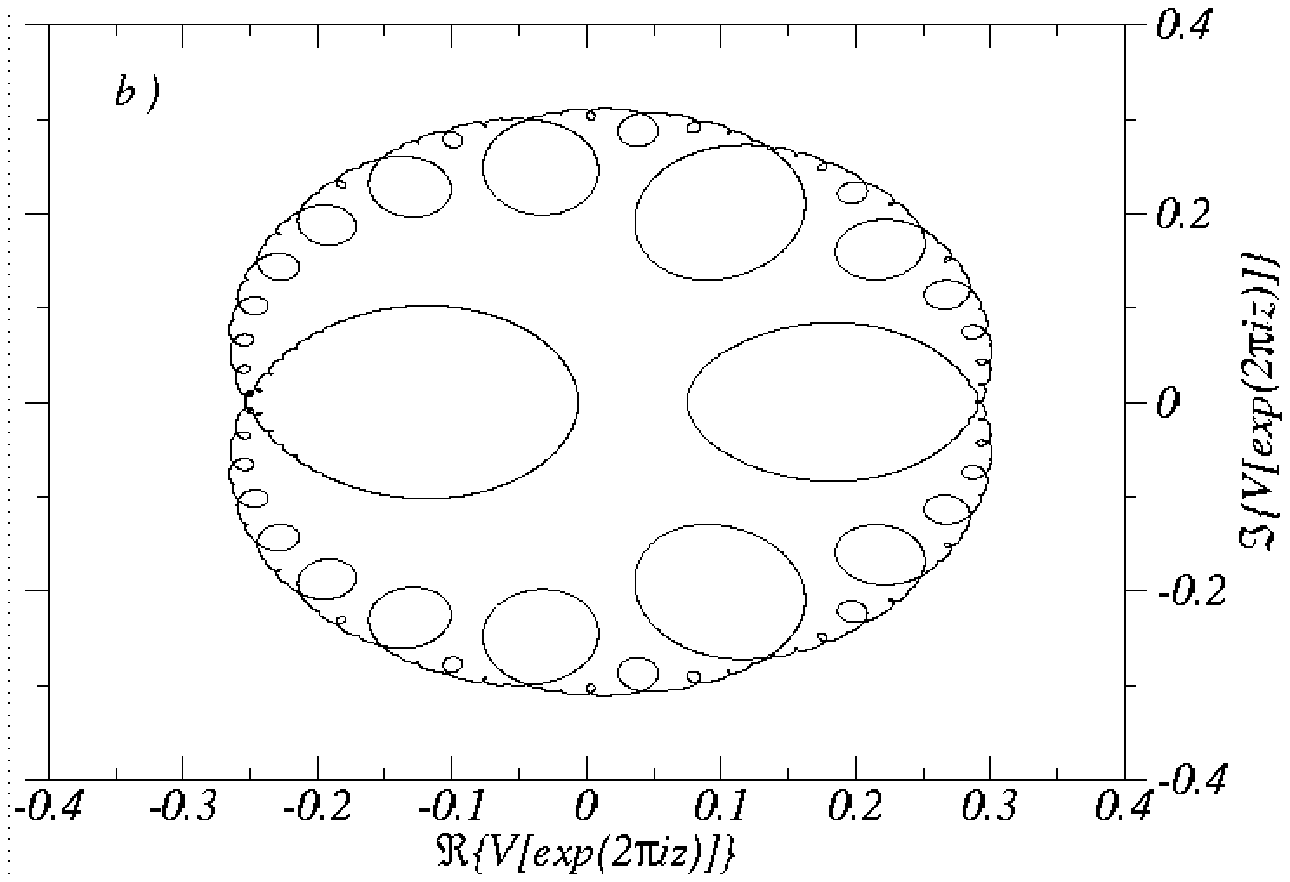}}
    }
    \end{center}
  \caption{Polar plot of $V\!\left( e^{2\pi i z}\right)$ for fixed values of $\Im z$.
  a)~$\Im z= 10^{-2}$ and b)~$\Im z= 10^{-3}$. $12000$ points
  uniformly distributed in $[0,1]$, $\epsilon_{V}=10^{-3}$.} 
  \label{fig:polarV}
  \end{figure}
\end{center}
\end{remark}

\section{The Littlewood--Paley Theory}
\label{sec:lpm}

The aim of this section is to introduce the basic ideas and results of the {\em Littlewood--Paley}
Theory, for a more complete discussion we refer to~\cite{Stein,Frazier} and also to~\cite{LlavePetrov} where
authors apply this Theory to study the regularity properties of the conjugating function for {\em critical circle maps}.
In \S~\ref{ssec:lpm} we will present the numerical implementation of this Theory to study the H\"older regularity of
the function $\mathcal{H}$ and the obtained estimate for the H\"older exponent.

The decay rate of the coefficients of a {\em Trigonometric series}, $\sum_{\Z}c_k e^{2\pi i kx}$, does not determine whether
this series is the {\em Fourier series} of some $L^p$ function if $p\neq 2$. More precisely given $f \in L^p$, $1\leq p <2$
and its Fourier series $\sum_{\Z}\hat{f}_k e^{2\pi i kx}$, then for ``almost every choice of signs $\pm 1$'', the series
$\sum_{\Z}(\pm 1)\hat{f}_k e^{2\pi i kx}$ is not the Fourier series of a $L^p$ function. This problem has been overcome by
Littlewood and Paley by ``grouping together'' trigonometric coefficients in {\em dyadic blocks}. Let
$A>1$, $\left(\mathcal{L}_0f\right)(x)=\hat{f}_0$ and, for $M\geq 1$, 
let $\left(\mathcal{L}_Mf\right)(x)=\sum_{A^{M-1}\leq |n| < A^{M}}\hat{f}_n e^{2\pi i n x}$, be the 
{\em dyadic partial sum} of $f$. Introducing the {\em Littlewood--Paley d--function}: 
$d(f)(x)=\left[ \sum_{M\geq 0} |\mathcal{L}_M f(x)|^2 \right]^{1/2}$, one can prove~\cite{LittlewoodPaley,Frazier} that
for all $1<p<+\infty$ there exist positive constants $A_p$ and $B_p$ such that:
\begin{equation*}
  A_p ||f||_p \leq ||d(f)||_p \leq B_p ||f||_p \, .
\end{equation*}

The Littlewood--Paley Theory is indeed more general, allowing to characterize other functional spaces by 
property of Fourier coefficients, for instance it applies~\cite{Frazier} to Sobolev spaces, Hardy spaces, 
H\"older spaces and Besov spaces. In the case of H\"older regularity one can easily realize that Fourier 
coefficients of an $\eta$--H\"older continuous 
function decay according to $\hat{f}_l =\mathcal{O}(|l|^{-\eta})$, the converse is not true but again the 
Littlewood--Paley Theory can characterize the H\"older regularity by the decay rate of the dyadic blocks.

An important tool in the Theory of Fourier series is the {\em Poisson Kernel}: $P_s(x)=\sum_{k\in \Z}s^{|k|}e^{2\pi i kx}$,
$s\in [0,1)$ and $x\in \T$. Let $(f*g)(x)=\int_0^1 f(\xi)g(x-\xi) \, d\xi$ be the convolution product for 
$1$--periodic functions. Then one can prove the following result (\cite{Stein} Lemma 5 
or \cite{Krantz} Theorem 15.6)
\begin{theorem}[Continuous Littlewood--Paley]
\label{thm:clp}
  Let $0 < \eta < r$, $r\in \N$ and $f$ be a continuous $1$--periodic function. Then $f$ is $\eta$--H\"older continuous
if and only if there exists $C>0$ such that for all $t>0$:
\begin{equation*}
\Big\lvert\Big\lvert\left(\frac{\partial}{\partial t}\right)^{r} 
\mathcal{P}_f(x,t)\Big\rvert\Big\rvert_{\infty} \leq C t^{\eta -r} \, ,
\end{equation*}
where $\mathcal{P}_f(x,t)=(P_{exp\left(-2\pi t\right)}*f)(x)$.
\end{theorem}
We remark that if the Theorem holds for some $r\in N$, then the same is true for any $r_1\in \N$, $r_1 > r$.
We call this Theorem {\em Continuous Littlewood--Paley} to distinguish it from the following result, which is more close
to the original idea of dyadic decomposition and we will call it {\em Discrete Littlewood--Paley} (see~\cite{Krantz} Theorem 5.9) 

\begin{theorem}[Discrete Littlewood-Paley]
\label{thm:dlp}
Let $\eta >0$ and let $f\in \mathcal{C}^0\left(\T\right)$. Then $f$
is $\eta$--H\"older continuous function if and only if for all $A>1$ there exists a
positive constant $C$ such that for all $M\in\N$ we have:
\begin{equation*}
||\mathcal{L}_Mf||_{\infty} \leq C A^{-\eta M} \, .
\end{equation*}
\end{theorem}
One usually take $A=2$, and so the name dyadic decomposition, but the result is independent of
the value of $A$. In the numerical implementation of this method we will use a value $A$ close to $1.25$
for computational reasons.

\section{Presentation of numerical results}
\label{sec:numrel}

This section collects our numerical results about the Marmi--Moussa--Yoccoz
conjecture that we recall here:
  \begin{quotation}
    {\em The analytic function, defined on $\H$: 
   $z \mapsto \mathcal{H}\left( z\right)=
  \log U \! (e^{2\pi i z}) -i \B(z)$, extends to a
  $1/2$--H\"older continuous function on the closure of $\H$.}
  \end{quotation}

\begin{center}
  \begin{figure}[ht]
    \begin{center}
    \mbox{\subfigure
    {\includegraphics[scale=0.3,angle=-90]{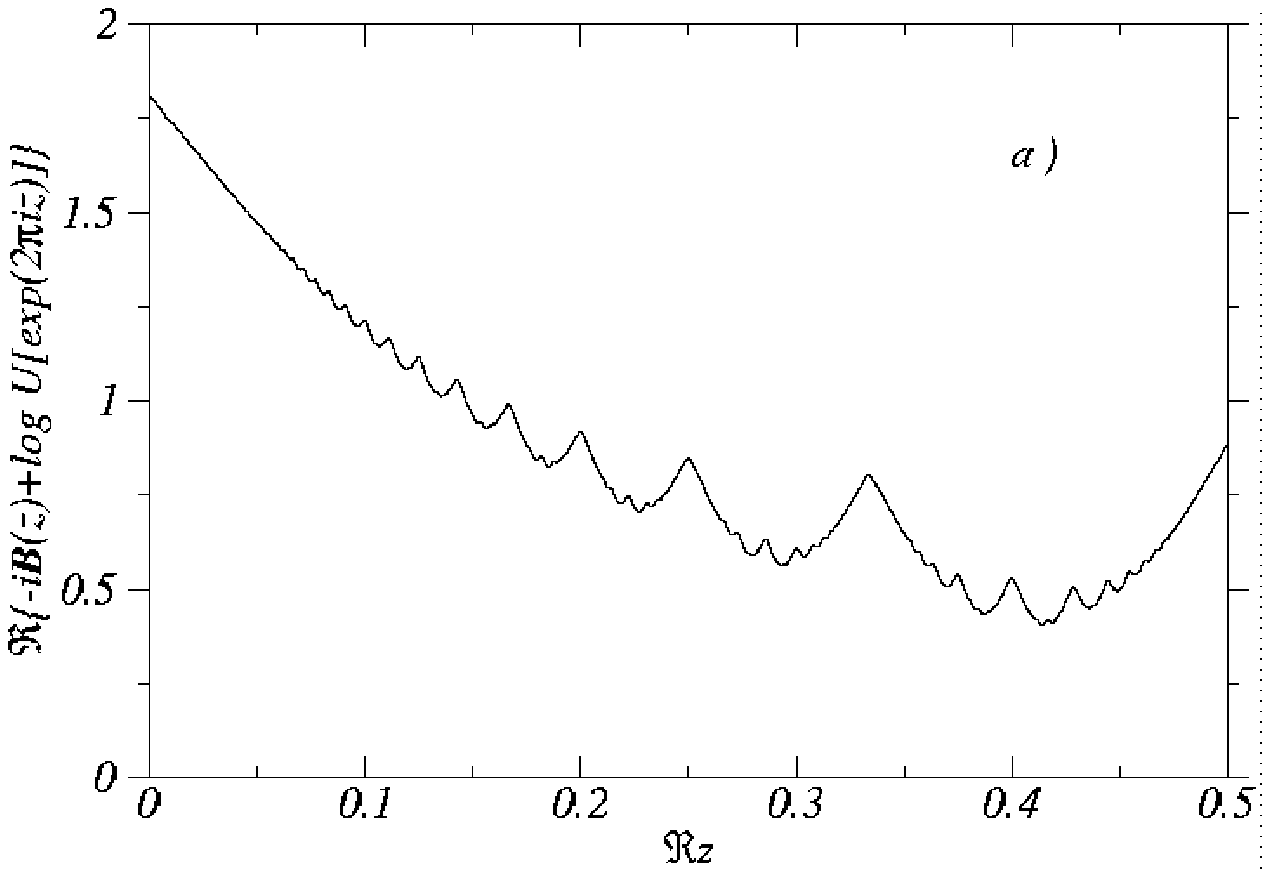}}\quad\quad
    \subfigure
    {\includegraphics[scale=0.3,angle=-90]{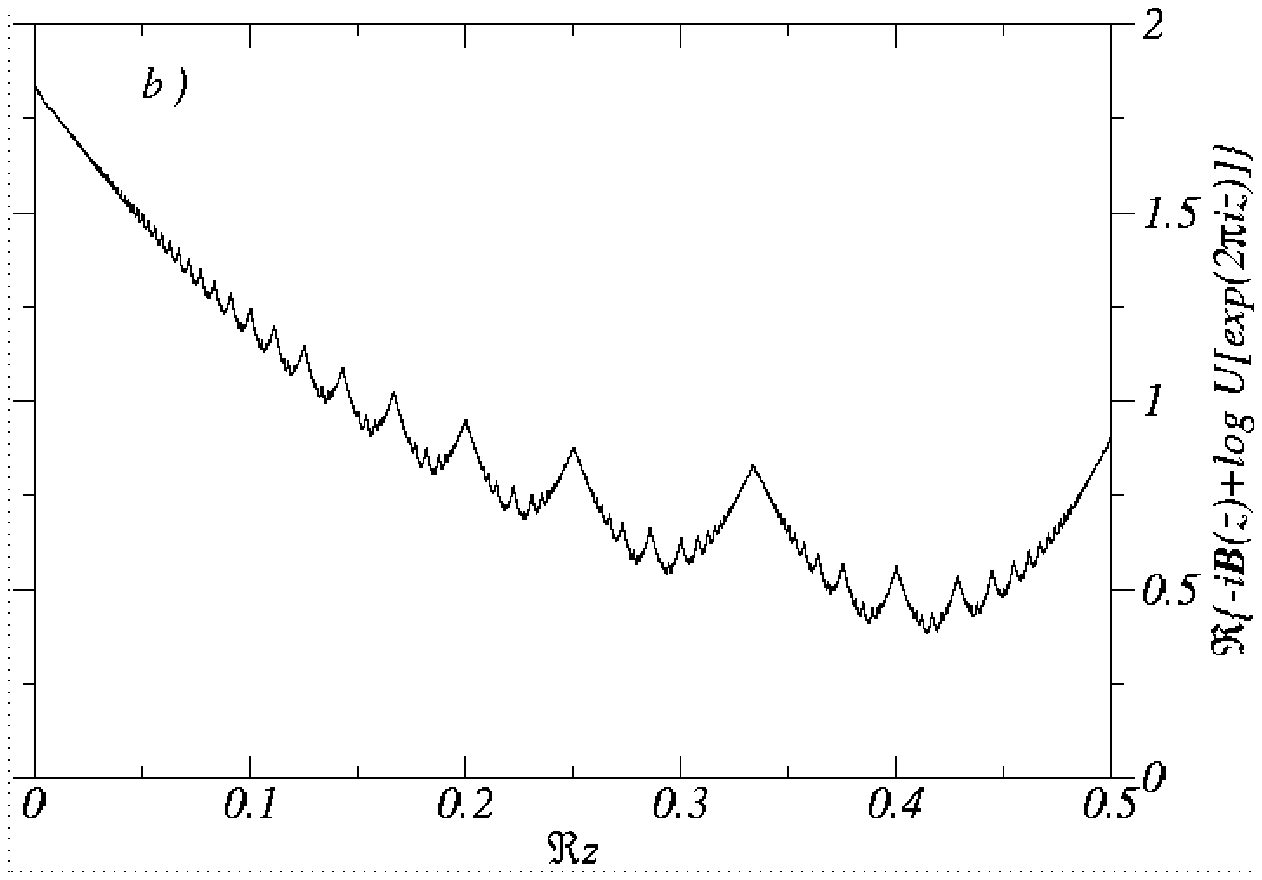}}
       }
    \mbox{\subfigure
    {\includegraphics[scale=0.3,angle=-90]{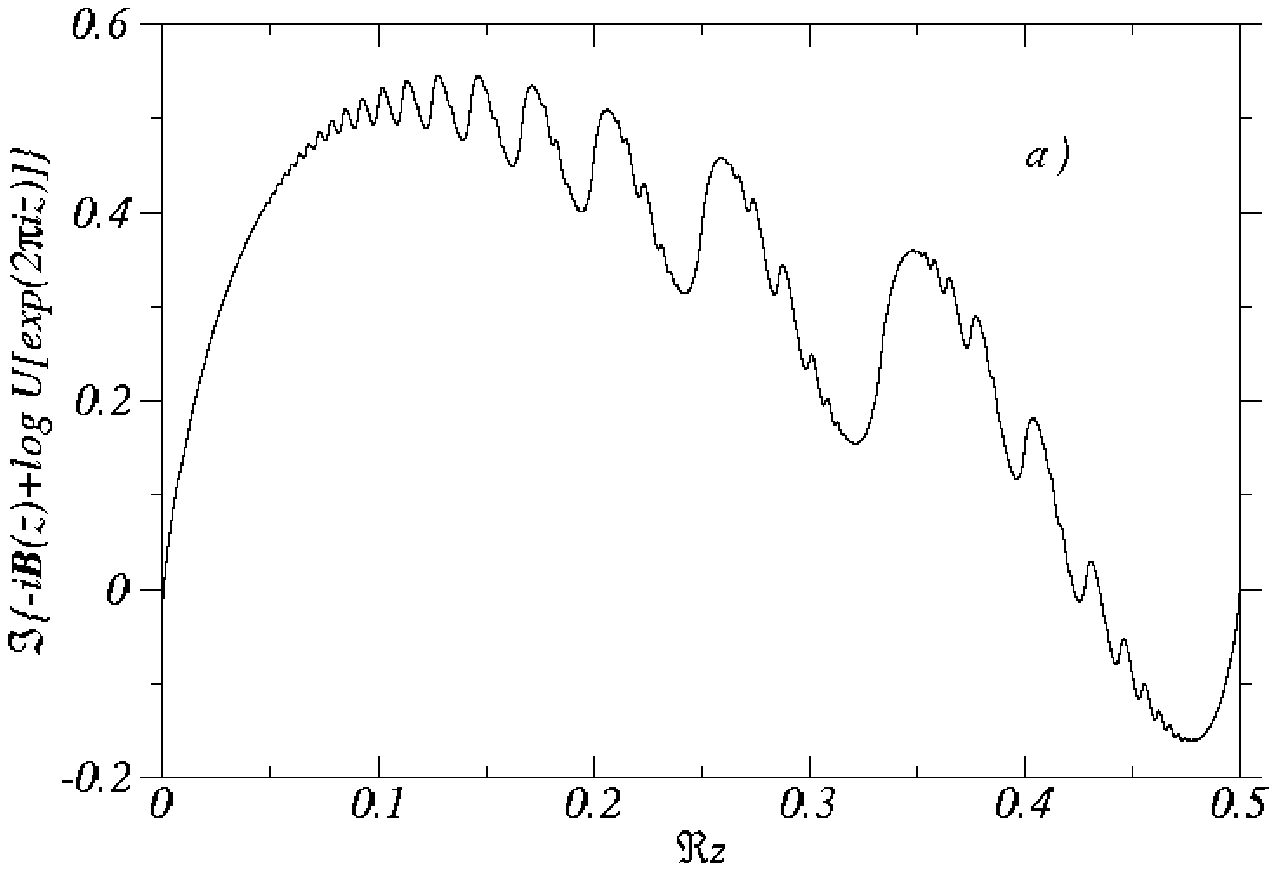}}\quad\quad
    \subfigure
     {\includegraphics[scale=0.3,angle=-90]{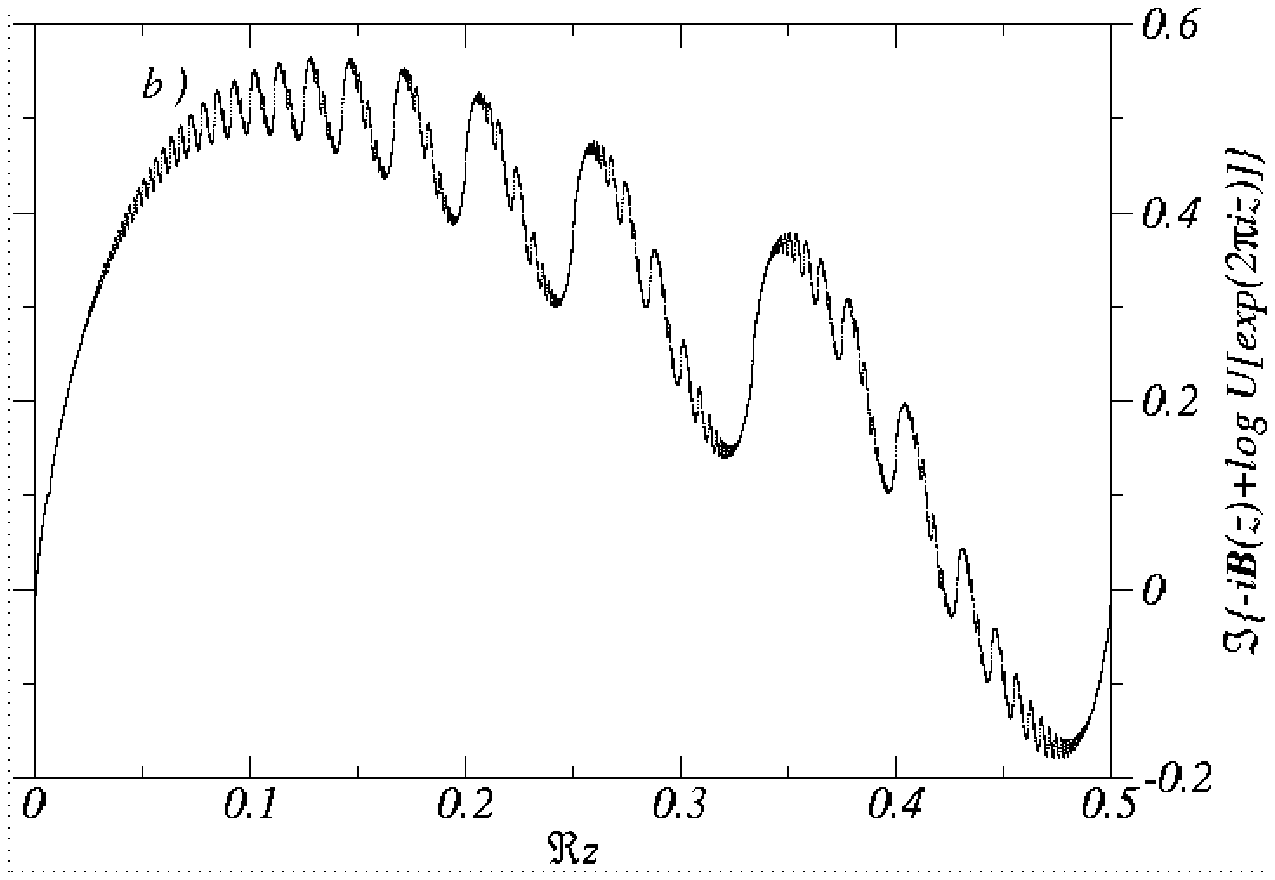}}
       }
    \end{center}
  \caption{Plot of $-i\B(z)+\log U\!\left( e^{2\pi i
  z}\right)$ vs $\Re z$ at fixed $\Im z$. On the top we show the real
  part whereas on the bottom the imaginary one. The first column a) is
  for~$\Im z= 10^{-3}$ and the second one b) for $\Im z= 10^{-4}$.
  Each plot has $10000$ points uniformly  distributed in
  $[0,1/2]$. $k_1=80$, $k_2=20$, $N_{max}=151$ and $\epsilon_U=10^{-3}$.} 
  \label{fig:immiBlogU}
  \end{figure}
\end{center}

Let us begin with some consideration concerning $\mathcal{H}$.
Remark~\ref{rem:fourcoeffB} and Remark~\ref{rem:parityfourcoeff} 
imply that $\mathcal{H}(z)$ has {\em even real part} for $\Re z \in
[0,1/2]$ and 
an {\em odd imaginary part} on the same domain. Moreover $\mathcal{H}(z)$ is
$1$--periodic and its Fourier series has only real coefficients, which
correspond to non--negative Fourier modes: 
$\mathcal{H}(z)=\sum_{l\geq 0}\hat{h}_le^{2\pi i lz}$.

\begin{center}
  \begin{figure}[ht]
    \begin{center}
    \mbox{\subfigure
    {\includegraphics[scale=0.3,angle=-90]{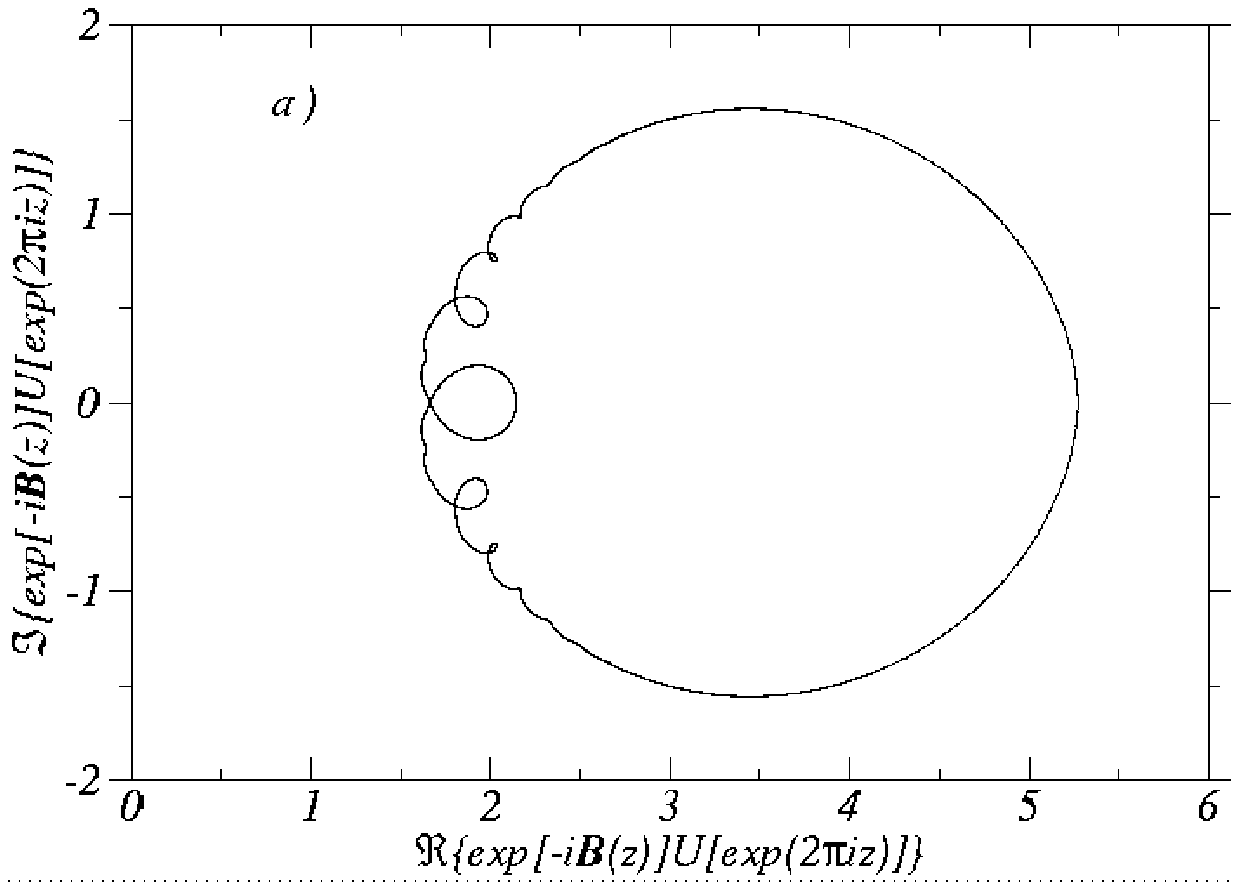}}\quad\quad
    \subfigure
    {\includegraphics[scale=0.3,angle=-90]{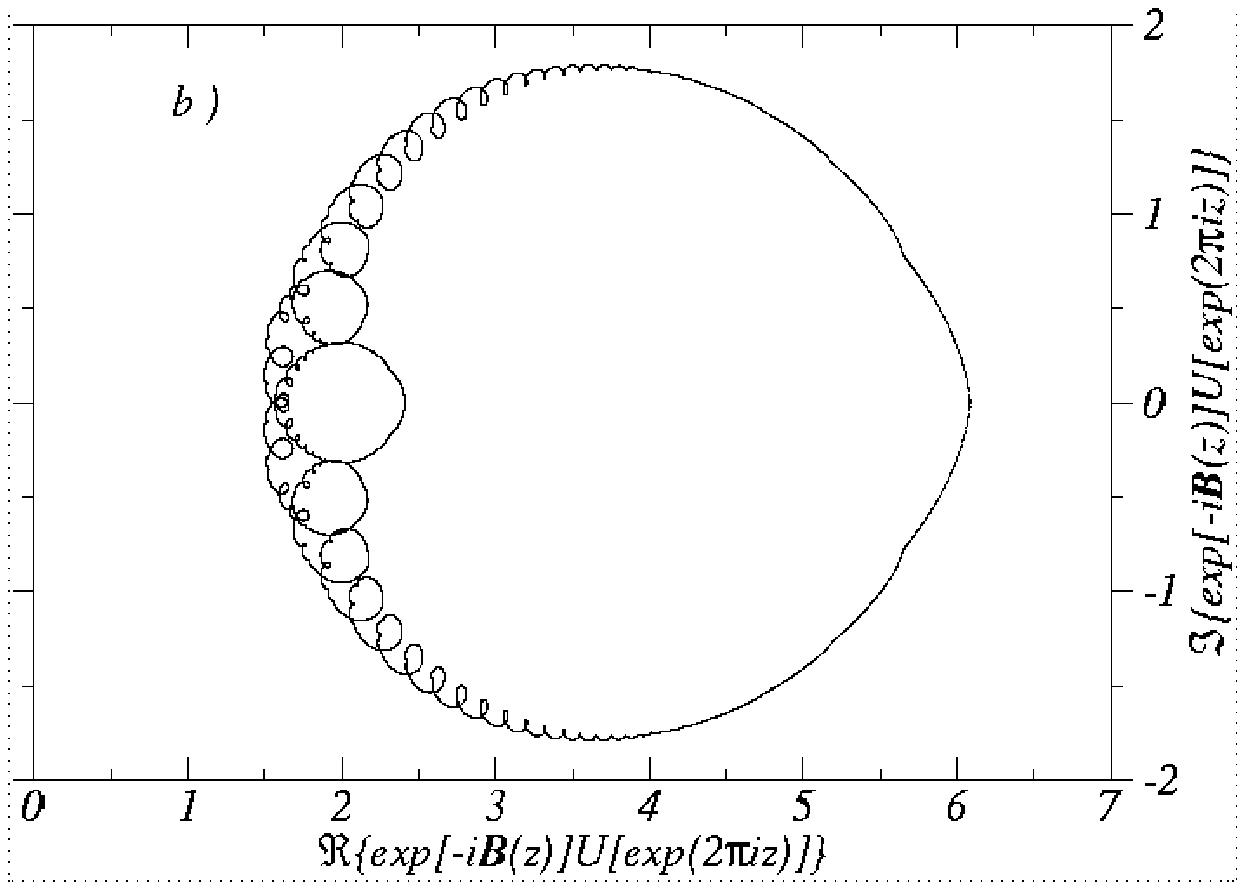}}
    }
    \end{center}
  \caption{Polar plot of $e^{-i\B\left(z\right)}U\!\left( e^{2\pi i
  z}\right)$ for fixed values of $\Im z$. 
  a)~$\Im z= 10^{-2}$ and b)~$\Im z= 10^{-3}$. $12000$ points
  uniformly distributed in $[0,1]$, $k_1=80$, $k_2=20$, $N_m=101$,
  $\epsilon_{U}=10^{-3}$.}  
  \label{fig:polaremiBU}
  \end{figure}
\end{center}

In Figure~\ref{fig:immiBlogU} we plot real and imaginary parts of
$\mathcal{H}(z)$ for some fixed small $\Im z$; remark that $\mathcal{H}$ still
has a ``structure'' but jumps of $\B$ and $-\log U$ seem to ``compensate''
to give a {\em continuous function}. The same fact holds for the
``bubbles'' of $U(e^{2\pi i z})$ and $e^{i \B(z)}$,
Figure~\ref{fig:polaremiBU} show some polar plots of
$e^{\mathcal{H}(z)}$ for fixed small $\Im z>0$, there are
still some ``bubbles'' but they are far from $(0,0)$. 

The H\"older continuity will be proved in the next paragraph, by giving
an estimate of the H\"older exponent applying the
 {\em Littlewood--Paley Theory} in the {\em Discrete}
and {\em Continuous} versions.

\subsection{Numerical Littlewood--Paley method.}
\label{ssec:lpm}

To numerically implement the Littlewood--Paley Theory we compute
from the numerical values of $\mathcal{H}$, a finite but large, number of
Fourier coefficients. Then to apply Theorem~\ref{thm:clp} we remark that
the convolution with the Poisson Kernel and the $r$--th derivative
have the form:
\begin{equation*}
  \left(\frac{\partial}{\partial t}\right)^{r} \left(P_{exp(-2\pi t)}*\mathcal{H}\right)(x)=
\sum_{l\geq 0} (-2\pi l)^re^{-2\pi tl}\hat{h}_l e^{2\pi i lx} \, ,
\end{equation*}
where we used the previous remark on the Fourier coefficients of $\mathcal{H}$.

We numerically compute $\Big\lvert\Big\lvert \left(\frac{\partial}{\partial t}\right)^{r} 
\left(P_{exp(-2\pi t)}*\mathcal{H}\right)\Big\rvert\Big\rvert_{\infty}$ for several small values of $t$
and some $r>1$, then applying a linear regression over the data:
\begin{equation}
\label{eq:clp}
  \log \Big\lvert\Big\lvert \left(\frac{\partial}{\partial t}\right)^{r} 
\left(P_{exp(-2\pi t)}*\mathcal{H}\right)\Big\rvert\Big\rvert_{\infty} = C_r^{\prime}-\beta_{CLP}(r) \log t \, ,
\end{equation}
we obtain a numerical value for $\eta^{(r)}_{CLP}=r-\beta_{CLP}(r)$.

From a numerical point of view the continuous version of the Littlewood--Paley method is better than the
discrete one, in fact the former has two parameters to vary $t$ and
$r$. We can vary $r$ to control
if the computed value of $\eta^{(r)}_{CLP}$ stays constant or not. Moreover we can compute the l.h.s. of~\eqref{eq:clp} 
for many values of $t$ and for each one all the known Fourier coefficients are 
used, whereas in the dyadic decomposition to ``small'' $M$ only ``few'' Fourier coefficients give their contribution and
only for ``large'' $M$ a large number of Fourier coefficients enter.

\begin{center}
  \begin{figure}[ht]
    \mbox{\subfigure
    {\includegraphics[scale=0.24]{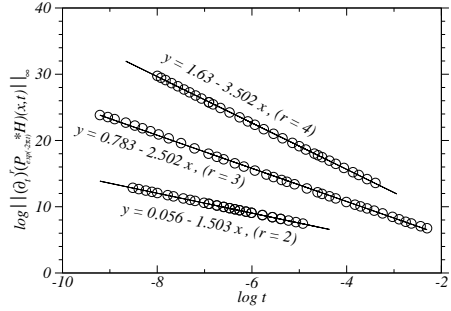}}
    }
  \caption{The function $\log t \mapsto 
\log \Big\lvert\Big\lvert\left( \frac{\partial}{\partial t}\right)^{r} 
 \mathcal{P}_{\mathcal{H}}(x,t) \Big\rvert\Big\rvert_{\infty}$,
for $r = 2$, $r=3$ and $r=4$. We also show the linear
regressions~\eqref{eq:clp}.}
  \label{fig:clp}
  \end{figure}
\end{center}

In Figure~\ref{fig:clp} we report data from~\eqref{eq:clp} and the
corresponding linear regression values~\footnote{In the Figure we
  decided to show only few points to have an ``intelligible picture'',
  but the linear regression are made using hundred of points.}. The
estimated values of $\eta$ obtained for different $r$ are:
$\eta_{CLP}^{(r=2)}=0.497 \pm 0.003$, $\eta_{CLP}^{(r=3)}=0.498 \pm 0.004$
and $\eta_{CLP}^{(r=4)}=0.498 \pm 0.003$ (errors are
standard deviation errors of linear regression). They agree in the
numerical precision and this gives a good indication of the validity
of the results. There is no reason to prefer one value to the other
and so we estimate $\eta_{CLP}=0.498 \pm 0.004$: the mean value of
the interval obtained by the union of the three intervals obtained for
$r=2,3,4$.

\begin{center}
  \begin{figure}[ht]
    \mbox{\subfigure
    {\includegraphics[scale=0.24]{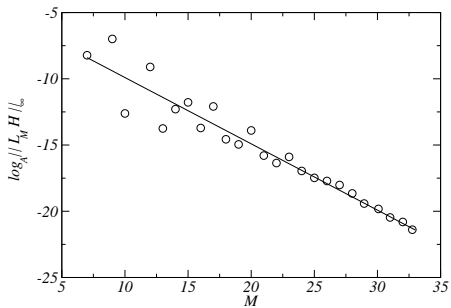}}
    }
  \caption{The function $M \mapsto
 \log_{A}||\mathcal{L}_M\mathcal{H}||_{\infty}$ and the linear fit 
$\log_{A}||\mathcal{L}_M\mathcal{H}||_{\infty}=C_{DLP}-\eta_{DLP}M$.}   
  \label{fig:dlp}
  \end{figure}
\end{center}

We also report the numerical results obtained using the discrete
Littlewood--Paley Theorem. 
We fix some $A>1$~\footnote{The exact value of $A$ is fixed in such a way we can take $M$ sufficiently large
to have a good asymptotic, even if we have a finite number of Fourier coefficients.} and from the computed 
Fourier coefficients of $\mathcal{H}$ we construct the dyadic partial sums for some large $M\in \N$. Then we use
a linear regression on the data:
\begin{equation}
\label{eq:dlp}
\log_A||\mathcal{L}_Mf||_{\infty} = C_{DLP} -\eta_{DLP} M \, ,
\end{equation}
to obtain the estimate value of the H\"older coefficient: $\eta_{DLP}=0.50\pm 0.03$ and $C_{DLP}=-4.90\pm 0.66$.
In Figure~\ref{fig:dlp} we report data from~\eqref{eq:dlp} and the linear regression applied on ``large $M$''.

\subsection{Conclusion}
\label{ssec:conc}

We conclude this paper by summarizing the obtained results. We
introduced the $1/2$--complex Bruno function and the Yoccoz function,
both with an algorithm to evaluate them numerically. We studied the
function $\mathcal{H}(z)=-i\B(z)+\log U\left( e^{2\pi i z}\right)$
defined on the upper Poincar\'e plane, and we conclude that it can be
extended to its closure, with a trace $\eta$--H\"older continuous.
Numerical results based on the Littlewood--Paley Theory give us the estimated
value for the H\"older exponent: $\eta_{CLP}=0.498$ with an error of $\pm 0.004$.
We can then conclude, with a good numerical evidence, that the
Marmi--Moussa--Yoccoz Conjecture should hold 
with the maximal exponent $1/2$.

\appendix

\section{Numerical considerations}
\label{sec:numcons}

The aim of this Appendix is to consider in details some technical parts of
 our numerical calculations. We will consider
 the role of the cut--off and their relations with the accuracy of the computations.
We will also compare the numerical properties of $\B$ with the analytical ones
proved in~\cite{MMYc}.

\subsection{Accuracy of the algorithm for $\B(z)$} 
\label{ssec:cutoff}

Let us recall the formula defining the $1/2$--complex Bruno function:
\begin{equation*}
\B(z)=\sum_{n\in \Z}\left[\sum_{g\in\hat{\mathcal{M}}}L_g\left(1+L_{\sigma}\right)\right] 
\varphi_{1/2}(z-n) \, ,
\end{equation*}
as already observed we need to introduce three cut--off to
compute it: $N_{max}$, $k_1$ and $k_2$. The first one determine the largest
Farey Series involved, namely only fractions $p/q$ s.t. $p/q\in [0,1)$,
 $(p,q)=1$ and $q\leq N_{max}$ will be considered to compute $\B$. The
other two cut--off: $k_1 \geq k_2 >0$, are introduced to truncate the sum
over $\Z$. Because the larger is $q$, the smaller is its contribution to
$\B$, to gain CPU times we decide to truncate the sum over $\Z$ at
$|n|\leq k_1$ if $q$ is ``small'', and to 
$|n|\leq k_2$ if $q$ is ``large''. Results showed in Section~\ref{sec:numrel} 
are obtained with $N_{max}=151$,
 $k_1=80$ and $k_2=20$. 

\begin{center}
  \begin{figure}[ht]
    \begin{center}
    \mbox{\subfigure
    {\includegraphics[scale=0.3,angle=-90]{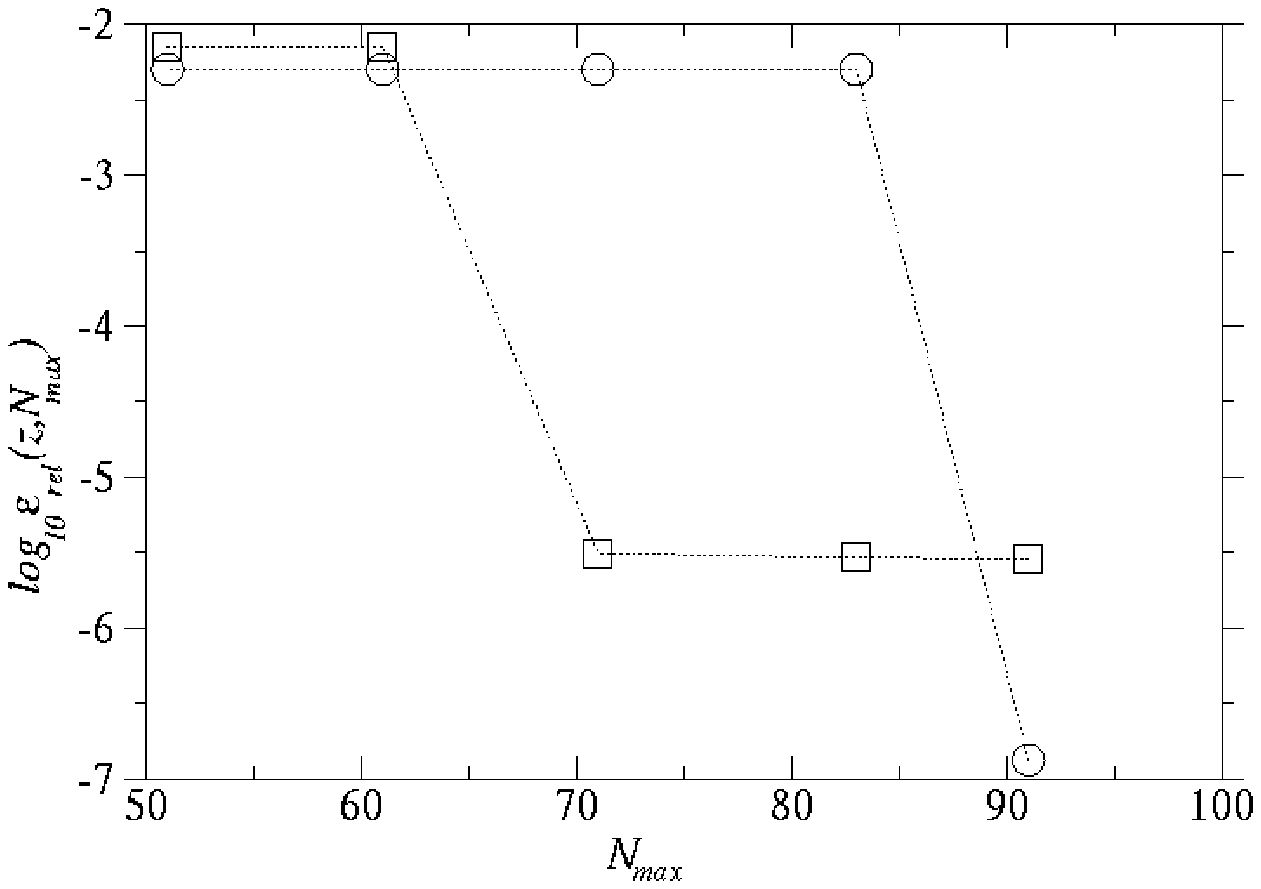}}\quad\quad
    \subfigure
    {\includegraphics[scale=0.3,angle=-90]{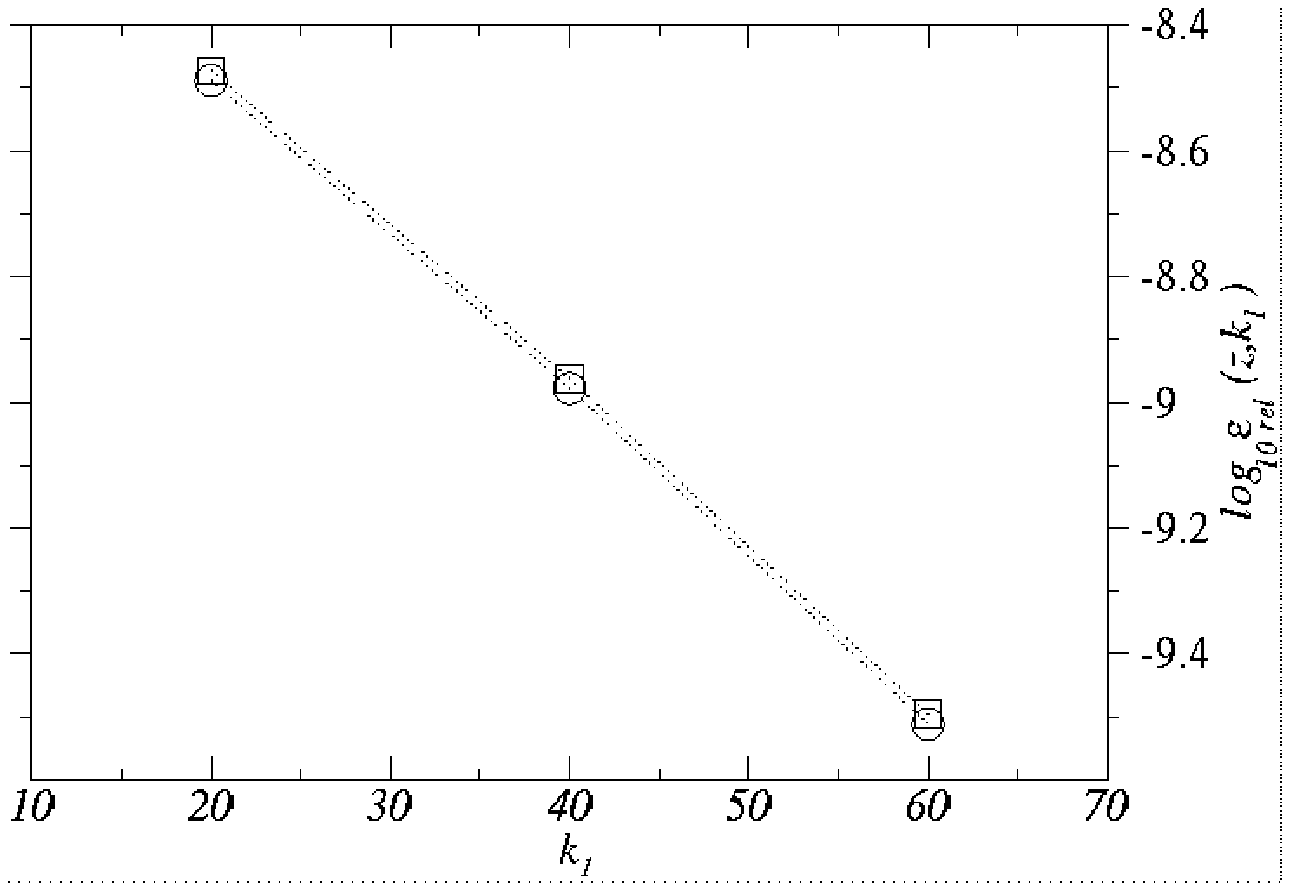}}}
    \end{center}
  \caption{Plot of $\log_{10}\epsilon_{rel}\left(z,M\right)$
  for some ``good'' $z$. On the left we show
  $\log_{10}\epsilon_{rel}\left(z,N_{max} \right)$ 
  whereas on the right $\log_{10}\epsilon_{rel}\left(z,k_1
  \right)$. Circles are for $z=2-\mathcal{G}+i 10^{-7}$ 
  and squares are for $z=\sqrt{2}-1+i 10^{-7}$.}
  \label{fig:cutoff}
  \end{figure}
\end{center}

In the rest of this paragraph we will study the
dependence of the computed Bruno function on the cut--off. Let us fix two
cut--off but one, call it generically $M$, and let us stress the
dependence of $\B$ on it, by setting $\B_M(z)$. We are then interested in 
studying the relative error: $\epsilon_{rel}\left(z,M \right)=
\lvert \B_M \left(z\right)-\B\left(z\right)\rvert/ \lvert \B \left(z\right)\rvert$,
where $z$ is fixed and $\B\left(z\right)$ is numerically computed with some
fixed large cut--off: $N_{max}=101$, $k_1=80$ and $k_2=20$. 
Or we can consider $\bar{\epsilon}_{rel}\left(M \right)$ the mean
value of $\epsilon_{rel}\left(z,M \right)$ for $\Re z \in [0,1/2]$ and
some fixed value of $\Im z >0$. In Table~\ref{tab:sup} we report values of 
$\log_{10} \bar{\epsilon}_{rel}\left(M \right)$, for $M\in \{
N_{max},k_1,k_2 \}$. Whereas in Figure~\ref{fig:cutoff}
 we show, $\log_{10}\epsilon_{rel}\left(z,M \right)$, 
for $M=N_{max}$ and $M=k_1$ and $z\in \{\sqrt{2}-1,2-\mathcal{G}
\}$. Clearly the larger are the cut--off the 
more accurate are the results, but recall that large cut--off implies large 
CPU times; in particular the CPU times increases almost linearly w.r.t. $k_1$ and $k_2$,
but quadratically w.r.t. $N_{max}$.

\begin{table}[ht]
\begin{tabular}[c]{|c|c||c|c||c|c|}
 \hline
 $N_{qmax}$ & $\log_{10} \bar{\epsilon}_{rel}\left(N_{max} \right)$ &
 $k_1$ & $\log_{10} \bar{\epsilon}_{rel}\left(k_1 \right)$ &
 $k_2$ & $\log_{10} \bar{\epsilon}_{rel}\left(k_2 \right)$ \\
 \hline\hline
 83 & -5.90 & 60 & -9.57 & 15 & -10.80 \\
 61 & -4.73 & 40 & -9.10 & 10 & -10.32 \\
 41 & -3.81 & 20 & -8.63 &    &        \\
 \hline
\end{tabular}
\label{tab:sup}
\caption{We report $\log_{10} \bar{\epsilon}_{rel}\left(N_{max} \right)$,
  $\log_{10} \bar{\epsilon}_{rel}\left(k_{1} \right)$ and $\log_{10}
  \bar{\epsilon}_{rel}\left(k_{2} \right)$.}  
\end{table}

\begin{center}
  \begin{figure}[ht]
    \begin{center}
    \mbox{\subfigure
    {\includegraphics[scale=0.3,angle=-90]{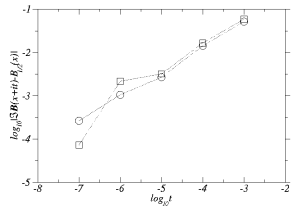}}}
    \end{center}
  \caption{Convergence of $\Im \B\left( x+it \right)$ to
  $B_{1/2}\left(x\right)$ when $t\rightarrow 0$ for some ``good
  irrational' $x$. Circles are for $x=2-\mathcal{G}$ whereas squares for
  $x=\sqrt{2}-1$. We plot $\log_{10}\lvert \Im \B\left( x+it
  \right)-B_{1/2}\left(x\right) \rvert$ versus 
  $\log_{10}t$. $N_{qmax}=101$, $k_1=80$ and $k_2=20$.}  
  \label{fig:limimag}
  \end{figure}
\end{center}

To have a full test of our algorithm we try to evaluate the limit, for
$\Im z \rightarrow 0$, of the computed $\B$ and compare it with the results 
proved in~\cite{MMYc}: section 5.2.9 page 816 and Theorem 5.19 page 827.
In particular we will be interested in studying the rate of convergence of
$\Im \B \left( x+it \right)$ to $B_{1/2}\left( x \right)$, for $t\rightarrow 0$ 
when $x$ is some ``good'' number (Figure~\ref{fig:limimag}), and 
the ``jump value'' of $\Re \B \left( p/q+it \right)$, when $t$ is
``small'' (Table A.2).

\begin{table}[ht]
\label{tab:limimag}
\begin{tabular}[c]{|c|c|}
 \hline
 $p/q$ & $\Delta \Re\B\left(p/q + it\right)-\pi/q$ \\
 \hline\hline
 $0/1$ & $1.1 \, 10^{-3}$ \\
 $1/2$ & $7.4 \, 10^{-4}$ \\
 $1/3$ & $3.6 \, 10^{-3}$ \\
 $1/4$ & $3.6 \, 10^{-3}$ \\
 $1/5$ & $3.8 \, 10^{-3}$ \\
 $2/5$ & $2.5 \, 10^{-3}$ \\
 \hline
\end{tabular}
  \caption{The jumps of $\Re\B\left(x+it\right)$ for rational $x$ and
  small $t$. The jump at $x=p/q$ is the numerical difference
  $|\Re\B\left(p/q+\delta+it\right)-\Re\B\left(p/q-\delta+it\right)|$,
  for $\delta$ small. We report the difference of the jump w.r.t the
  expected value for $x\in \{ 0/1,1/2,1/3,1/4,1/5,2/5 \}$ and $t=10^{-7}$.} 
\end{table}

\end{document}